\documentclass[11pt]{article}
\usepackage{amsmath}
\usepackage{bm}
\usepackage[numbers]{natbib}
\usepackage[colorlinks,citecolor=blue,linkcolor=blue,urlcolor=blue]{hyperref}
\usepackage{latexsym, epsfig, amssymb, amsmath, amsthm, graphicx, mathrsfs, booktabs}
\usepackage{rotating, tabularx, setspace,paralist}
\usepackage{color}
\usepackage{graphicx}
\usepackage[OT1]{fontenc}
\usepackage{lscape}
\usepackage{multirow,multicol}
\usepackage{anysize}
\renewcommand{\baselinestretch}{1.2}
\marginsize{1.1in}{1.1in}{0.55in}{0.8in} %

\makeatletter
\def\singlespace{\def\baselinestretch{1}\@normalsize}

\newtheorem{assumption}{Condition}
\newtheorem{lemma}{Lemma}
\newtheorem{proposition}{Proposition}
\newtheorem{theorem}{Theorem}
\newtheorem{definition}{Definition}
\renewcommand{\theequation}{
\arabic{equation}%
}

\makeatother

\newcommand{\ba}{\mbox{\bf a}}
\newcommand{\bb}{\mbox{\bf b}}
\newcommand{\be}{\mbox{\bf e}}

\newcommand{\bg}{\mbox{\bf g}}

\newcommand{\bs}{\mbox{\bf s}}

\newcommand{\bx}{\mbox{\bf x}}
\newcommand{\by}{\mbox{\bf y}}
\newcommand{\bz}{\mbox{\bf z}}
\newcommand{\bA}{\mbox{\bf A}}
\newcommand{\bB}{\mbox{\bf B}}
\newcommand{\bC}{\mbox{\bf C}}

\newcommand{\bG}{\mbox{\bf G}}
\newcommand{\bH}{\mbox{\bf H}}
\newcommand{\bI}{\mbox{\bf I}}

\newcommand{\bU}{\mbox{\bf U}}
\newcommand{\bM}{\mbox{\bf M}}

\newcommand{\bW}{\mbox{\bf W}}
\newcommand{\bX}{\mbox{\bf X}}

\newcommand{\bZ}{\mbox{\bf Z}}
\newcommand{\bone}{\mbox{\bf 1}}
\newcommand{\bzero}{\mbox{\bf 0}}
\newcommand{\bveps}{\mbox{\boldmath $\varepsilon$}}

\newcommand{\bbeta}{\mbox{\boldmath $\beta$}}
\newcommand{\brho}{\mbox{\boldmath $\rho$}}

\newcommand{\bdelta}{\mbox{\boldmath $\delta$}}
\newcommand{\btheta}{\mbox{\boldmath $\theta$}}

\newcommand{\bleta}{\mbox{\boldmath $\eta$}}

\newcommand{\bSig}{\mbox{\boldmath $\Sigma$}}
\newcommand{\bOmg}{\mbox{\boldmath $\Omega$}}
\newcommand{\tbx}{\widetilde \bx}
\newcommand{\tbX}{\widetilde \bX}
\newcommand{\hbX}{\widehat \bX}

\newcommand{\wt}{\widetilde}
\newcommand{\mb}{\mathbb}
\newcommand{\wh}{\widehat}

\newcommand{\hbbeta}{\widehat\bbeta}

\newcommand{\hbeta}{\widehat{\beta}}

\newcommand{\calA}{{\cal A}}

\newcommand{\sbOmg}{\mbox{\scriptsize\boldmath $\Omega$}}

\newcommand{\cov}{\mathrm{cov}}

\newcommand{\Sig}{\mathbf{\Sigma}}
\newcommand{\veps}{\varepsilon}
\newcommand{\tr}{\mathrm{tr}}
\newcommand{\diag}{\mathrm{diag}}

\newcommand{\supp}{\mathrm{supp}}

\newcommand{\argmin}{\mbox{argmin}}

\def\t{^T}

\newcommand\independent{\protect\mathpalette{\protect\independenT}{\perp}}
\def\independenT#1#2{\mathrel{\setbox0\hbox{$#1#2$}%
\copy0\kern-\wd0\mkern4mu\box0}}

\begin{document}

\title{RANK: Large-Scale Inference with Graphical Nonlinear Knockoffs
\thanks{Yingying Fan is Associate Professor, Data Sciences and Operations Department, Marshall School of Business, University of Southern California, Los Angeles, CA 90089 (E-mail: \textit{fanyingy@marshall.usc.edu}). Emre Demirkaya is Ph.D. Candidate, Department of Mathematics, University of Southern California, Los Angeles, CA 90089 (E-mail: \textit{demirkay@usc.edu}). Gaorong Li was Postdoctoral Scholar at USC during this work and is Professor, Beijing Institute for Scientific and Engineering Computing, Beijing University of Technology, Beijing, China 100124 (E-mail: \textit{ligaorong@gmail.com}). Jinchi Lv is McAlister Associate Professor in Business Administration, Data Sciences and Operations Department, Marshall School of Business, University of Southern California, Los Angeles, CA 90089 (E-mail: \textit{jinchilv@marshall.usc.edu}). This work was supported by NSF CAREER Award DMS-1150318 and a grant from the Simons Foundation.}
\date{\today}
\author{Yingying Fan$^1$, Emre Demirkaya$^1$, Gaorong Li$^2$ and Jinchi Lv$^1$
\medskip\\
University of Southern California$^1$ and Beijing University of Technology$^2$
\\
} %
}

\maketitle

\begin{abstract}
Power and reproducibility are key to enabling refined scientific discoveries in contemporary big data applications with general high-dimensional nonlinear models. In this paper, we provide theoretical foundations on the power and robustness for the model-free knockoffs procedure introduced recently in Cand\`{e}s, Fan, Janson and Lv (2016) in high-dimensional setting when the covariate distribution is characterized by Gaussian graphical model. We establish that under mild regularity conditions, the power of the oracle knockoffs procedure with known covariate distribution in high-dimensional linear models is asymptotically one as sample size goes to infinity.  When moving away from the ideal case, we suggest the modified model-free knockoffs method called graphical nonlinear knockoffs (RANK) to accommodate the unknown covariate distribution. We provide theoretical justifications on the robustness of our modified procedure by showing that the false discovery rate (FDR) is asymptotically controlled at the target level and the power is asymptotically one with the estimated covariate distribution. To the best of our knowledge, this is the first formal theoretical result on the power for the knockoffs procedure. Simulation results demonstrate that compared to existing approaches, our method performs competitively in both FDR control and power. A real data set is analyzed to further assess the performance of the suggested knockoffs procedure.
\end{abstract}

\textit{Running title}: RANK

\textit{Key words}: Power; Reproducibility; Big data; High-dimensional nonlinear models; Robustness; Large-scale inference and FDR; Graphical nonlinear knockoffs

\section{Introduction} \label{sec:Intro}

Feature selection with big data is of fundamental importance to many contemporary applications from different disciplines of social sciences, health sciences, and engineering \cite{HastieTibshiraniFriedman2009, FanLv2010, BuhlmannvandeGeer2011}. Over the past two decades, various feature selection methods, theory, and algorithms have been extensively developed and investigated for a wide spectrum of flexible models ranging from parametric to semiparametric and nonparametric linking a high-dimensional covariate vector $\bx = (X_1, \cdots, X_p)\t$ of $p$ features $X_j$'s to a response $Y$ of interest, where the dimensionality $p$ can be large compared to the available sample size $n$ or even greatly exceed $n$. The success of feature selection for enhanced prediction in practice can be attributed to the reduction of noise accumulation associated with high-dimensional data through dimensionality reduction. In particular, most existing studies have focused on the power perspective of feature selection procedures such as the sure screening property, model selection consistency, oracle property, and oracle inequalities. When the model is correctly specified, researchers and practitioners often would like to know whether the estimated model involving a subset of the $p$ covariates enjoys reproducibility in that the fraction of noise features in the discovered model is controlled. Yet such a practical issue of reproducibility is largely less well understood for the settings of general high-dimensional nonlinear models. Moreover, it is no longer clear whether the power of feature selection procedures can be retained when one intends to ensure the reproducibility.

Indeed, the issues of power and reproducibility are key to enabling refined scientific discoveries in big data applications utilizing general high-dimensional nonlinear models. To characterize the reproducibility of statistical inference, the seminal paper of \cite{BenjaminiHochberg1995} introduced an elegant concept of false discovery rate (FDR) which is defined as the expectation of the fraction of false discoveries among all the discoveries, and proposed a popularly used Benjamini--Hochberg 
procedure for FDR control by resorting to the p-values for large-scale multiple testing returned by some statistical estimation and testing procedure. There is a huge literature on FDR control for large-scale inference and various generalizations and extensions of the original FDR procedure were developed and investigated for different settings and applications \cite{BenjaminiYekutieli2001, EfronTibshirani2002, Storey2002, StoreyTaylorSiegmund2004, AbramovichBenjaminiDonohoJohnstone2006, Efron2007JASA, Efron2007AOS, FanHallYao2007, Wu2008, ClarkeHall2009, HallWang2010, FanFan2011, MengSunZhangWaterman2011, ZhangLiu2011, FanHanGu2012, LiuShao2014, SuCandes2016}. Most of existing work either assumes a specific functional form such as linearity
on the dependence structure of response $Y$ on covariates $X_j$'s, or relies on the p-values for evaluating the significance of covariates $X_j$'s. Yet in high-dimensional settings, we often do not have such luxury since response $Y$ could depend on covariate vector $\bx$ through very complicated forms and even when $Y$ and $\bx$ have simple dependence structure, high dimensionality of covariates can render classical p-value calculation procedures no longer justified or simply invalid \cite{Huber1973, FanDemirkayaLv2017, SurChenCandes2017}. These intrinsic challenges can make the p-value based methods difficult to apply or even fail \cite{CandesFanJansonLv2016}.

To accommodate arbitrary dependence structure of $Y$ on $\bx$ and bypass the need of calculating accurate p-values for covariate significance, \cite{CandesFanJansonLv2016} recently introduced the model-free knockoffs framework for FDR control in general high-dimensional nonlinear models. Their work was inspired by and builds upon the ingenious development of the knockoff filter in \cite{BarberCandes2015}, which provides effective FDR control in the setting of Gaussian linear model with dimensionality $p$ no larger than sample size $n$. The knockoff filter was later extended in \cite{BarberCandes2016} to high-dimensional linear model using the ideas of data splitting and feature screening. The salient idea of \cite{BarberCandes2015} is to construct the so-called ``knockoff" variables which mimic the dependence structure of the original covariates but are independent of response $Y$ conditional on the original covariates. These knockoff variables can be used as control variables. By comparing the regression outcomes for original variables with those for control variables, the relevant set of variables can be  identified more accurately and thus the FDR can be better controlled. The model-free knockoffs framework introduced in \cite{CandesFanJansonLv2016} greatly expands the applicability of the original knockoff filter in that the response $Y$ and covariates $\bx$ can have arbitrarily complicated dependence structure and the dimensionality $p$ can be arbitrarily large compared to sample size $n$. It was theoretically justified in \cite{CandesFanJansonLv2016} that the model-free knockoffs procedure controls FDR exactly in finite samples of arbitrary dimensions. However, one important assumption in their theoretical development is that the joint distribution of covariates $\bx$ should be known. Moreover, formal power analysis of the knockoffs framework is still lacking even for the setting of Gaussian linear model.

Despite the importance of known covariate distribution in their theoretical development, \cite{CandesFanJansonLv2016} empirically explored the scenario of unknown covariate distribution for the specific setting of generalized linear model (GLM) \cite{McCullaghNelder1989} with Gaussian design matrix and discovered that the estimation error of the covariate distribution can have  negligible effect on FDR control. Yet there exist no formal theoretical justifications on the robustness of the model-free knockoffs method and it is also unclear to what extent such robustness can hold beyond the GLM setting. To address these fundamental challenges, our paper intends as the first attempt to provide theoretical foundations on the power and robustness for the model-free knockoffs framework. Specifically, the major innovations of the paper are twofold. First, we will provide theoretical support on the robustness of the model-free knockoffs procedure with unknown covariate distribution in general high-dimensional nonlinear models. Second, we will formally investigate the power of the knockoffs framework in high-dimensional linear models with both known and unknown covariate distribution.

More specifically, in the ideal case of known covariate distribution we prove that the model-free knockoffs procedure in \cite{CandesFanJansonLv2016} has asymptotic power one under mild regularity conditions in high-dimensional linear models. When moving away from the ideal scenario,
to accommodate the difficulty caused by unknown covairate distribution we suggest the modified model-free knockoffs method called graphical nonlinear knockoffs (RANK). The modified knockoffs procedure exploits the data splitting idea, where the first half of the sample is used to estimate the unknown covariate distribution and reduce the model size, and the second half of the sample is employed to \textit{globally} construct the knockoff variables and apply the knockoffs procedure. We establish that the modified knockoffs procedure asymptotically controls the FDR regardless of whether the reduced model contains the true model or not. Such feature makes our work intrinsically different from that in \cite{BarberCandes2016} requiring the sure screening property \cite{FanLv2008} of the reduced model; see Section \ref{sec: mod-knockoffs} for more detailed discussions on the differences.
In our theoretical analysis of FDR, we still allow for arbitrary dependence structure of response $Y$ on covariates $\bx$ and assume that the joint distribution of $\bx$ is characterized by Gaussian graphical model with unknown precision matrix \cite{Lauritzen1996}. In the specific case of high-dimensional linear models with unknown covariate distribution, we also provide robustness analysis on the power of our modified procedure.

The rest of the paper is organized as follows. Section \ref{Sec2} reviews the model-free knockoffs framework and provides theoretical justifications on its power in high-dimensional linear models. We introduce the modified model-free knockoffs procedure RANK and investigate its robustness on both FDR control and power with respect to the estimation of unknown covariate distribution in Section \ref{sec3}. Section \ref{sec4} presents several simulation examples of both linear and nonlinear models to verify our theoretical results. We demonstrate the performance of our procedure on a real data set in Section \ref{sec5}. Section \ref{sec6} discusses some implications and extensions of our work. The proofs of main results are relegated to the Appendix. Additional technical details are provided in the Supplementary Material.

\section{Power analysis for oracle model-free knockoffs} \label{Sec2}


Suppose we have a sample $(\bx_i, Y_i)_{i = 1}^n$ of $n$ independent and identically distributed (i.i.d.) observations from the population $(\bx, Y)$, where dimensionality $p$ of covariate vector $\bx = (X_1, \cdots, X_p)\t$ can greatly exceed available sample size $n$. To ensure model identifiability, it is common to assume that only a small fraction of $p$ covariates $X_j$'s are truly relevant to response $Y$. To be more precise, \cite{CandesFanJansonLv2016} defined the set of irrelevant features $\mathcal{S}_1$ as that consisting of $X_j$'s such that $X_j$ is independent of $Y$ conditional on all remaining $p - 1$ covariates $X_k$'s with $k \neq j$, and thus the set of truly relevant features $\mathcal{S}_0$ is given naturally by $\mathcal{S}_1^c$, the complement of set $\mathcal{S}_1$. Features in sets $\mathcal{S}_0$ and $\mathcal{S}_0^c = \mathcal{S}_1$ are also referred to as important and noise features, respectively.

We aim at accurately identifying these truly relevant features in set $\mathcal{S}_0$ while keeping the false discovery rate (FDR) \cite{BenjaminiHochberg1995} under control. The FDR for a feature selection procedure is defined as
\begin{align}\label{eq-fdr}
& {\rm FDR}=\mathbb{E}[{\rm FDP}] \quad \text{ with } {\rm FDP}=\dfrac{|\widehat{\mathcal{S}} \cap \mathcal{S}_0^c |}{|\widehat{\mathcal{S}}|},
\end{align}
where $\widehat{\mathcal{S}}$ denotes the sparse model returned by the feature selection procedure, $|\cdot|$ stands for the cardinality of a set, and the convention $0/0 = 0$ is used in the definition of the false discovery proportion (FDP) which is the fraction of noise features in the discovered model. Here feature selection procedure can be any favorite sparse modeling method by the choice of the user.

\subsection{Review of model-free knockoffs framework} \label{subsec2.1}

Our suggested graphical nonlinear knockoffs procedure in Section \ref{sec3} falls in the general framework of model-free knockoffs introduced in \cite{CandesFanJansonLv2016}, which we briefly review in this section. The key ingredient of model-free knockoffs framework is the construction of the so-called model-free knockoff variables that are defined as follows.

\begin{definition}[\cite{CandesFanJansonLv2016}] \label{def1}
Model-free knockoffs for the family of random variables $\bx = (X_1,\cdots, X_p)^T$ is a new family of random variables $\widetilde{\bx} = (\widetilde{X}_1, \cdots, \widetilde{X}_p)^T$ that satisfies two properties: (1) $(\bx^T, \widetilde{\bx}^T)_{\emph{swap}(\mathcal{S})} {\overset{d}{=}}(\bx^T, \widetilde{\bx}^T)$ for any subset $\mathcal{S} \subset\{1,\cdots, p\}$,
	where $\emph{swap}(\mathcal{S})$ means swapping components $X_j$ and $\widetilde{X}_j$ for each $j\in \mathcal{S}$ and $\overset{d}{=}$ denotes equal in distribution, and (2) $\widetilde{\bx} \independent Y | \bx$.
\end{definition}

We see from Definition \ref{def1} that model-free knockoff variables $\widetilde{X}_j$'s mimic the probabilistic dependency structure among the original features $X_j$'s and are independent of response $Y$ given $X_j$'s. When the covariate distribution is characterized by Gaussian graphical model \cite{Lauritzen1996}, that is,
\begin{equation} \label{neweq004}
\bx\sim N(\bzero, \bOmg_0^{-1})
\end{equation}
with $p \times p$ precision matrix $\bOmg_0$ encoding the graphical structure of the conditional dependency among the covariates $X_j$'s, we can construct the $p$-variate model-free knockoff random variable $\tbx$ characterized in Definition \ref{def1} as
\begin{equation}\label{equ-cond-ko}
\tbx|\bx \sim N\Big(\bx - \diag\{\bs\} \bOmg_0 \bx, 2\diag\{\bs\} - \diag\{\bs\} \bOmg_0 \diag\{\bs\}\Big),
\end{equation}
where $\bs$ is a $p$-dimensional vector with nonnegative components chosen in a suitable way. In fact, in view of (\ref{neweq004}) and (\ref{equ-cond-ko}) it is easy to show that the original features and model-free knockoff variables have the following joint distribution
\begin{align} \label{eq007}
\left(\begin{array}{c}
\bx \\
\tbx
\end{array}\right) \sim N\left(\left(\begin{array}{c}
\bzero \\
\bzero
\end{array}\right), \left(\begin{array}{cc}
\bSig_0 & \bSig_0 - \diag\{\bs\} \\
\bSig_0 - \diag\{\bs\} & \bSig_0
\end{array} \right) \right)
\end{align}
with $\bSig_0 = \bOmg_0^{-1}$ the covariance matrix of covariates $\bx$. Intuitively, larger components of $\bs$ means that the constructed knockoff variables deviate further from the original features, resulting in higher power in distinguishing them. The $p$-dimensional vector $\bs$ in (\ref{equ-cond-ko}) should be chosen in a way such that $\bSig_0 - 2^{-1} \diag\{\bs\}$ is positive definite, and can be  selected using the methods in \cite{CandesFanJansonLv2016}. We will treat it as a nuisance parameter throughout our theoretical analysis.

With the constructed knockoff variables $\tbx$, the knockoffs inference framework proceeds as follows. We select important variables by resorting to the knockoff statistics $W_j = f_j(Z_j, \wt{Z}_j)$ defined for each $1 \leq j \leq p$, where $Z_j$ and $\wt{Z}_j$ represent feature importance measures for $j$th covariate $X_j$ and its knockoff counterpart $\widetilde{X}_j$, respectively, and $f_j(\cdot,\cdot)$ is an antisymmetric function satisfying $f_j(z_j,\tilde{z}_j) = -f_j(\tilde{z}_j, z_j)$. For example, in linear regression models one can choose $Z_j$ and $\wt{Z}_j$ as the Lasso \cite{Tibshirani1996} regression coefficients of $X_j$ and $\widetilde{X}_j$, respectively, and a valid knockoff statistic is $W_j = f_j(z_j, \tilde z_j) = |z_j| - |\tilde z_j|$. There are also many other options for defining the feature importance measures. Observe that all model-free knockoff variables $\widetilde{X}_j$'s are just noise features by the second property in Definition \ref{def1}. Thus intuitively, a large positive value of knockoff statistic $W_j$ indicates that $j$th covariate $X_j$ is important, while a small magnitude of $W_j$ usually corresponds to noise features.

The final step of the knockoffs inference framework is to sort $|W_j|$'s from high to low and select features whose $W_j$'s are at or above some threshold $T$, which results in the discovered model
\begin{equation} \label{neweq005}
	\widehat{\mathcal{S}} =	\widehat{\mathcal{S}}( T)=\{1\leq j \leq p: W_j\geq  T\}.
\end{equation}
Following \cite{BarberCandes2015} and \cite{CandesFanJansonLv2016}, one can choose the threshold $T$ in the following two ways
\begin{eqnarray}
T&=&\min\Bigg\{t\in\mathcal{W}: \dfrac{|\{j:  W_j\le-t\}|}{|\{j:  W_j\geqslant t\}|}\le q\Bigg\},\label{T-knock}\\
T_{+}&=&\min\Bigg\{t\in\mathcal{W}: \dfrac{1+|\{j:  W_j\le-t\}|}{|\{j:  W_j\geqslant t\}|}\le q\Bigg\},\label{T-threshold+}
\end{eqnarray}
where $\mathcal{W}=\{|W_j|: 1 \leq j \leq p\} \setminus \{0\}$ is the set of unique nonzero values attained by $|W_j|$'s and $q \in (0, 1)$ is the desired FDR level specified by the user. The procedures using threshold $T$ in (\ref{T-knock}) and threshold $T_{+}$ in (\ref{T-threshold+}) are referred to as knockoffs and knockoffs$_+$ methods, respectively. It was proved in \cite{CandesFanJansonLv2016} that model-free knockoffs procedure controls a modified FDR that replaces $|\widehat{\mathcal{S}}|$ in the denominator by $q^{-1} + |\widehat{\mathcal{S}}|$ in (\ref{eq-fdr}), and model-free knockoffs$_+$ procedure achieves exact FDR control in finite samples regardless of dimensionality $p$ and dependence structure of response $Y$ on covariates $\bx$. The major assumption needed in \cite{CandesFanJansonLv2016} is that the distribution of covariates $\bx$ is known. Throughout the paper, we implicitly use the threshold $T_{+}$ defined in (\ref{T-threshold+}) for FDR control in the knockoffs inference framework but still write it as $T$ for notational simplicity.

\subsection{Power analysis in linear models} \label{sec: oracle-power}

Although the knockoffs procedures were proved rigorously to have controlled FDR in \cite{BarberCandes2015, BarberCandes2016, CandesFanJansonLv2016}, their power advantages over popularly used approaches have been demonstrated only numerically therein. In fact, formal power analysis for the knockoffs framework is still lacking even in simple model settings such as linear regression. We aim to fill in this gap as a first attempt and provide theoretical foundations on the power analysis for model-free knockoffs framework. In this section, we will focus on the \textit{oracle} model-free knockoffs procedure for the ideal case when the true precision matrix $\bOmg_0$ for the covariate distribution in (\ref{neweq004}) is known, which is the setting assumed in \cite{CandesFanJansonLv2016}. The robustness analysis for the case of unknown precision matrix $\bOmg_0$ will be undertaken in Section \ref{sec3}.

Since the power analysis for the knockoffs framework is nontrivial and challenging, we content ourselves on the setting of high-dimensional linear models for the technical analysis on power. The linear regression model assumes that
\begin{align}\label{eq001}
  \by = \bX\bbeta_0 + \bveps,
  \end{align}
where $\by=(Y_1,\cdots, Y_n)^T$ is an $n$-dimensional response vector, $\bX = (\bx_1, \cdots, \bx_n)^T$ is an $n \times p$ design matrix consisting of $p$ covariates $X_j$'s, $\bbeta_0 = (\beta_{0, 1}, \cdots, \beta_{0, p})\t$ is a $p$-dimensional true regression coefficient vector, and $\bveps = (\veps_1,\cdots, \veps_n)^T$ is an $n$-dimensional error vector independent of $\bX$. As mentioned before, the true model $\mathcal{S}_0 = \supp(\bbeta_0)$ which is the support of $\bbeta_0$ is assumed to be sparse with size $s=|\mathcal{S}_0|$, and the $n$ rows of design matrix $\bX$ are i.i.d. observations generated from Gaussian graphical model (\ref{neweq004}). Without loss of generality, all the diagonal entries of covariance matrix $\bSig_0$ are assumed to be ones.

As discussed in Section \ref{subsec2.1}, there are many choices of the feature selection procedure up to the user for producing the feature importance measures $Z_j$ and $\wt{Z}_j$ for covariates $X_j$ and knockoff variables $\widetilde{X}_j$, respectively, and there are also different ways to construct the knockoff statistics $W_j$. For the illustration purpose, we adopt the Lasso coefficient difference (LCD) as the knockoff statistics in our power analysis. The specific choice of LCD for knockoff statistics was proposed and recommended in \cite{CandesFanJansonLv2016}, in which it was demonstrated empirically to outperform some other choices in terms of power. The LCD is formally defined as
  \begin{eqnarray}\label{WH-statistics0}
  {W}_j= |\widehat{\beta}_j(\lambda)|-|\widehat{\beta}_{p+j}(\lambda)|,
  \end{eqnarray}
  where $\widehat{\beta}_j(\lambda)$ and $\widehat{\beta}_{p+j}(\lambda)$ denote the $j$th and $(p+j)$th components, respectively, of the Lasso \cite{Tibshirani1996} regression coefficient vector
  \begin{equation}\label{SLasso-solution0}
  \widehat{\bbeta}(\lambda) = {\argmin}_{\bb\in\mathbb{R}^{2p}}\left\{(2n)^{-1}\big\|\by-[\bX, {\widetilde{\bX}}]\bb\big\|_2^2+\lambda\|\bb\|_1\right\} 
  \end{equation}
  with   $\lambda\geq 0$  the regularization parameter, $\tbX = (\tbx_1,\cdots,\tbx_n)^T$ an $n \times p$ matrix whose $n$ rows are independent random vectors of model-free knockoff variables generated from (\ref{equ-cond-ko}), and $\|\cdot\|_r$ for $r \geq 0$ the $L_r$-norm of a vector.
To simplify the technical analysis, we assume that with asymptotic probability one, there are no ties in the magnitude of nonzero $W_j$'s and no ties in the magnitude of nonzero components of Lasso solution in \eqref{SLasso-solution0}, which is a mild condition in light of the continuity of the underlying distributions.

To facilitate the power analysis, we impose some basic regularity conditions.

  \begin{assumption}\label{cond4}
  	The components of $\bveps$ are i.i.d. with sub-Gaussian distribution.
  \end{assumption}

  \begin{assumption}\label{cond3}
  	It holds that $\{n/(\log p)\}^{1/2} \min\limits_{j\in \mathcal{S}_0}|\beta_{0,j}|\rightarrow \infty$ as $n$ increases.
  \end{assumption}

\begin{assumption}\label{cond:Shat0}
	There exists some constant $c \in (2(qs)^{-1},1)$  such that with asymptotic probability one, $|\widehat{\mathcal{S}}| \geq cs$ for $\widehat{\mathcal{S}}$ given in (\ref{neweq005}). 
\end{assumption}

Condition \ref{cond4} can be relaxed to heavier-tailed distributions at the cost of slower convergence rates as long as similar concentration inequalities used in the proofs continue to hold.  Condition \ref{cond3} is assumed to ensure that the Lasso solution $\hbbeta(\lambda)$ does not miss a great portion of important features in $\mathcal{S}_0$. This is necessary since the knockoffs procedure under investigation builds upon the Lasso solution and thus its power is naturally upper bounded by that of Lasso. To see this, recall the well-known oracle inequality for Lasso \cite{Bickel2009, BuhlmannvandeGeer2011} that with asymptotic probability one,
$
\|\hbbeta(\lambda) - \bbeta_0\|_2 = O(s^{1/2} \lambda)
$
for $\lambda$ chosen in the order of $\{(\log p)/n\}^{1/2}$. Then Condition \ref{cond3} entails that for some $\kappa_n \rightarrow \infty$,
$
O(s\lambda^2) = \|\hbbeta(\lambda) - \bbeta_0\|_2^2  \geq \sum_{j \in \wh{\mathcal S}_{\text{L}}^c\cap \mathcal{S}_0}\beta_{0,j}^2 \geq n^{-1} (\log p)\kappa_n^2 |\wh{\mathcal S}_{\text{L}}^c\cap \mathcal{S}_0|
$
with $\wh{\mathcal S}_{\text{L}} = \supp\{\hbbeta(\lambda)\}$.  Thus the number of important features missed by Lasso $|\wh{\mathcal S}_{\text{L}}^c\cap \mathcal{S}_0|$ is upper bounded by $O(s \kappa_n^{-2})$ with asymptotic probability one. This guarantees that the power of Lasso is lowered bounded by $1-O(\kappa_n^{-2})$; that is, Lasso has asymptotic power one.  We will show in Theorem \ref{theo1} that there is almost no power loss when applying model-free knockoffs procedure.

Condition \ref{cond:Shat0} imposes a lower bound on the size of the sparse model selected by the knockoffs procedure. Since it is not straightforward to check, we provide a sufficient condition that is more intuitive in Lemma \ref{lem: cond-Shat0} below, which shows that Condition \ref{cond:Shat0} can hold as long as there exist enough strong signals in the model. We acknowledge that Lemma \ref{lem: cond-Shat0} may not be a necessary condition for Condition \ref{cond:Shat0}.

\begin{lemma}\label{lem: cond-Shat0}
Assume that Condition \ref{cond4} holds and there exists some constant $c \in (2(qs)^{-1},1)$ such that $|\mathcal{S}_2|\geq cs$ with $\mathcal{S}_2 = \{j: |\beta_{0,j}|\gg [sn^{-1}(\log p)]^{1/2}\}$. Then Condition \ref{cond:Shat0} holds.
\end{lemma}

We are now ready to characterize the statistical power of the knockoffs procedure in high-dimensional linear model \eqref{eq001}. Formally speaking, the power of a feature selection procedure is defined as
\begin{equation}\label{power-def1}
{\rm Power}(\wh{\mathcal{S}})=\mb{E}\Big[\dfrac{|\widehat{\mathcal S} \cap \mathcal{S}_0|}{|\mathcal{S}_0|}\Big],
\end{equation}
where $\widehat{\mathcal{S}}$ denotes the discovered sparse model returned by the feature selection procedure.

\begin{theorem}\label{theo1} Assume that Condition \ref{cond4}--\ref{cond:Shat0} hold, all the eigenvalues of $\bOmg_0$ are bounded away from $0$ and $\infty$, the smallest eigenvalue of $2\diag(\bs) - \diag(\bs)\bOmg_0\diag(\bs)$ is positive and bounded away from $0$, and $\lambda= C_1\{(\log p)/n\}^{1/2}$ with $C_1 > 0$ some constant. Then the oracle model-free knockoffs procedure satisfies that with asymptotic probability one, 	
$|\widehat{\mathcal S} \cap \mathcal{S}_0| / |\mathcal{S}_0| \geq1 - O(\kappa_n^{-1})$
for some $\kappa_n\rightarrow\infty$, and ${\rm Power}(\wh{\mathcal{S}}) \rightarrow 1$ as $n \rightarrow \infty$.
\end{theorem}

Theorem \ref{theo1} reveals that the oracle model-free knockoffs procedure in \cite{CandesFanJansonLv2016} knowing the true precision matrix $\bOmg_0$ for the covariate distribution can indeed have asymptotic power one under some mild regularity conditions. This shows that for the ideal case, model-free knockoffs procedure can enjoy appealing FDR control and power properties simultaneously.

\section{Robustness of graphical nonlinear knockoffs} \label{sec3}
When moving away from the ideal scenario considered in Section \ref{Sec2}, a natural question is whether both properties of FDR control and power can continue to hold with no access to the knowledge of true covariate distribution. To gain insights into such a question, we now turn to investigating the robustness of model-free knockoffs framework. Hereafter we assume that the true precision matrix $\bOmg_0$ for the covariate distribution in (\ref{neweq004}) is \textit{unknown}. We will begin with the FDR analysis and then move on to the power analysis.

\subsection{Modified model-free knockoffs} \label{sec: mod-knockoffs}
We would like to emphasize that the linear model assumption is no longer needed here and arbitrary dependence structure of response $\by$ on covariates $\bx$ is allowed. As mentioned in Introduction, to overcome the difficulty caused by unknown precision matrix $\bOmg_0$ we modify the model-free knockoffs procedure described in Section \ref{subsec2.1} and suggest the method of graphical nonlinear knockoffs (RANK).

To ease the presentation, we first introduce some notation.
For each given $p \times p$ symmetric positive definite matrix $\bOmg$, denote by $\bC^{\sbOmg} = \bI_p- \diag\{\bs\}\bOmg$ and $\bB^{\sbOmg}= \big(2\diag\{\bs\} - \diag\{\bs\}\bOmg\diag\{\bs\}\big)^{1/2}$ the square root matrix. We define $n \times p$ matrix $\widetilde{\bX}^{\sbOmg} = (\widetilde{\bx}_1^{\sbOmg},\cdots, \widetilde{\bx}_n^{\sbOmg})^T$ by independently generating $\widetilde{\bx}_i^{\sbOmg}$ from the conditional distribution
\begin{equation}\label{eq: xtilde-distr}
\widetilde\bx_i^{\sbOmg}|\bx_i \sim N\Big(\bC^{\sbOmg}\bx_i, (\bB^{\sbOmg})^2\Big),
\end{equation}
where $\bX = (\bx_1, \cdots, \bx_n)^T$ is the original $n \times p$ design matrix generated from Gaussian graphical model (\ref{neweq004}). It is easy to show that the $(2p)$-variate random vectors $(\bx_i^T, (\widetilde{\bx}_i^{\sbOmg})^T)^T$ are i.i.d. with Gaussian distribution of mean $\bzero$ and covariance matrix given by $\cov(\bx_i)  = \bSig_0$, $\cov(\bx_i, \tbx_i^{\sbOmg}) = \bSig_0\bC^{\sbOmg}$, and $\cov(\tbx_i^{\sbOmg})  = (\bB^{\sbOmg})^2 + \bC^{\sbOmg}\bSig_0(\bC^{\sbOmg})^T$.

Our modified knockoffs method RANK exploits the idea of data splitting, in which one half of the sample is used to estimate unknown precision matrix $\bOmg_0$ and reduce the model dimensionality, and the other half of the sample is employed to construct the knockoff variables and implement the knockoffs inference procedure, with the steps detailed below.
\begin{itemize}
\item \textit{Step 1}. Randomly split the data $(\bX, \by)$ into two folds $(\bX^{(k)}, \by^{(k)})$ with $1 \leq k \leq 2$ each of sample size $n/2$.
	
\item \textit{Step 2}. Use the first fold of data $(\bX^{(1)}, \by^{(1)})$ to obtain an estimate $\wh{\bOmg}$ of the precision matrix and a reduced model with support $\wt{\mathcal{S}}$.
	
\item \textit{Step 3}. With estimated precision matrix $\wh{\bOmg}$ from Step 2, construct an $(n/2) \times p$ knockoffs matrix $\hbX$ using  $\bX^{(2)}$ with rows independently generated from (\ref{eq: xtilde-distr}); that is,
	$
	\hbX = \bX^{(2)}(\bC^{\wh{\sbOmg}})^{T} + \bZ (\bB^{\wh{\sbOmg}})^2
	$
	with $\bZ$ an $(n/2)\times p$ matrix with i.i.d. $N(0, 1)$ components.

\item \textit{Step 4}. Construct knockoff statistics $W_j$'s using only data on support $\wt{\mathcal{S}}$, that is,
	$W_j = W_j(\by^{(2)}, \bX_{\wt{\mathcal{S}}}^{(2)}, \hbX_{\wt{\mathcal{S}}})$ for $j\in \wt{\mathcal{S}}$
	  and $W_j = 0$ for $j\in \wt{\mathcal{S}}^c$. Then apply knockoffs inference procedure to $W_j$'s to obtain final set of features $\wh{\mathcal{S}}$.
\end{itemize}
Here for any matrix $\bA$ and subset $\mathcal{S}\subset \{1,\cdots, p\}$, the compact notation $\bA_{\mathcal{S}}$ stands for the submatrix of $\bA$ consisting of columns in set $\mathcal{S}$.

As discussed in Section \ref{subsec2.1}, the model-free knockoffs framework utilizes sparse regression procedures such as the Lasso. For this reason, even in the original model-free knockoffs procedure the knockoff statistics $W_j$'s (see, e.g., \eqref{WH-statistics0}) take nonzero values only over a much smaller model than the full model. This observation motivates us to estimate such a smaller model using the first half of the sample in Step 2 of our modified procedure.  When implementing this modified procedure, we limit ourselves to sparse models $\wt{\mathcal{S}}$ with size bounded by some positive integer $K_n$ that diverges with $n$; see, for example,  \cite{FanLv2013, Lv2013} for detailed discussions and justifications on similar consideration of sparse models. In addition to sparse regression procedures, feature screening methods such as \cite{FanLv2008, FanFan2008} can also be used to obtain the reduced model $\wt{\mathcal S}$.

The above modified knockoffs method differs from the original model-free knockoffs procedure \cite{CandesFanJansonLv2016} in that we use an independent sample to obtain the estimated precision matrix $\wh{\bOmg}$ and reduced model $\wt{\mathcal{S}}$. In particular, the independence between estimates $(\wh{\bOmg}, \wt{\mathcal{S}})$ and data $(\bX^{(2)}, \by^{(2)})$ plays an important role in our theoretical analysis for the robustness of the knockoffs procedure. In fact, the idea of data splitting has been popularly used in the literature for various purposes \cite{FanSamworthWu2009, FanGuoHao2012, ShahSamworth2013, BarberCandes2016}. Although the work of \cite{BarberCandes2016} has the closest connection to ours, there are several key differences between these two methods. Specifically, \cite{BarberCandes2016} considered high-dimensional linear model with fixed design, where the data is split into two portions with the first portion used for feature screening and the second portion employed for applying the original knockoff filter in \cite{BarberCandes2015} on the reduced model. To ensure FDR control, it was required in \cite{BarberCandes2016} that the feature screening method should enjoy the sure screening property \cite{FanLv2008}, that is, the reduced model after the screening step contains the true model $\mathcal{S}_0$ with asymptotic probability one. In contrast, one major advantage of our method is that the asymptotic FDR control can be achieved without requiring the sure screening property; see Theorem \ref{thm1} in Section \ref{sec: 2.2} for more details.  Such major distinction is rooted on the difference in constructing knockoff variables; that is, we construct model-free knockoff variables \textit{globally} in Step 3 above, whereas \cite{BarberCandes2016} constructed knockoff variables \textit{locally} on the reduced model. Another major difference is that our method works with random design and does not need any assumption on how response $\by$ depends upon covariates $\bx$, while the method in \cite{BarberCandes2016} requires the linear model assumption and cannot be extended to nonlinear models.

\subsection{Robustness of FDR control for graphical nonlinear knockoffs} \label{sec: 2.2}

We begin with investigating the robustness of FDR control for the modified model-free knockoffs procedure RANK. To simplify the notation, we rewrite $(\bX^{(2)}, \by^{(2)})$ as $(\bX, \by)$ with sample size $n$ whenever there is no confusion, where $n$ now represents half of the original sample size. For each given $p \times p$ symmetric positive definite matrix $\bOmg$, an $n \times p$ knockoffs matrix $\wt{\bX}^{\sbOmg} = (\widetilde{\bx}_1^{\sbOmg},\cdots, \widetilde{\bx}_n^{\sbOmg})^T$ can be constructed with $n$ rows independently generated according to (\ref{eq: xtilde-distr}) and the modified knockoffs procedure proceeds with a given reduced model $\mathcal{S}$. Then the FDP and FDR functions in (\ref{eq-fdr}) can be rewritten as
\begin{equation} \label{neweq006}
{\rm FDR}_n(\bOmg, \mathcal{S}) = {\mb E}[{\rm FDP}_n(\by, \bX_{\mathcal{S}}, \tbX_\mathcal{S}^{\sbOmg})],
\end{equation}
where the subscript $n$ is used to emphasize the dependence of FDP and FDR functions on sample size.
It is easy to check that the knockoffs procedure based on $(\by, \bX_{\mathcal{S}}, \tbX_\mathcal{S}^{\sbOmg_0})$ satisfies all the conditions  in \cite{CandesFanJansonLv2016} for FDR control for any reduced model $\mathcal{S}$ that is independent of $\bX$ and $\tbX^{\sbOmg_0}$, which ensures that ${\rm FDR}_n(\bOmg_0, \mathcal{S})$ can be controlled at the target level $q$. To study the robustness of our modified knockoffs procedure, we will make a connection between functions ${\rm FDR}_n(\bOmg, \mathcal{S})$ and  ${\rm FDR}_n(\bOmg_0, \mathcal{S})$.

To ease the presentation, denote by $\tbX_0 = \tbX^{\sbOmg_0}$ the oracle knockoffs matrix with $\bOmg = \bOmg_0$, $\bC_0 = \bC^{\sbOmg_0}$, and $\bB_0 = \bB^{\sbOmg_0}$. The following proposition establishes a formal characterization of the FDR as a function of the precision matrix $\bOmg$ used in generating the knockoff variables and the reduced model $\mathcal{S}$.

\begin{proposition}\label{prop1}
For any given symmetric positive definite $\bOmg\in \mathbb{R}^{p\times p}$ and $\mathcal{S}\subset\{1,\cdots, p\}$, it holds that
	\begin{equation}
	{\rm FDR}_n(\bOmg, \mathcal{S})	= \mb{E}\Big[g_{n}\big(\bX_{\rm aug}^\mathcal{S}\bH^{\sbOmg}\big)\Big], \label{eq: gn-def}
	\end{equation}
	where $ \bX_{\rm aug}^\mathcal{S} = [\bX, \tbX_{0,\mathcal{S}}] \in \mb{R}^{n\times (p+|\mathcal{S}|)}$, function $g_n(\cdot)$ is some conditional expectation of the FDP function whose functional form is free of $\bOmg$ and $\mathcal{S}$, and
	\begin{align*}
	\bH^{\sbOmg} = \left(\begin{array}{cc}
		\bI_{p} & \bC_{\mathcal S}^{\sbOmg}- \bC_{0,\mathcal{S}}(\bB_{0,\mathcal{S}}^T\bB_{0,\mathcal{S}})^{-1/2}\big((\bB^{\sbOmg}_{\mathcal{S}})^T\bB^{\sbOmg}_{\mathcal{S}}\big)^{1/2} \\
		\bzero & (\bB_{0,\mathcal{S}}^T\bB_{0,\mathcal{S}})^{-1/2}\big((\bB^{\sbOmg}_{\mathcal{S}})^T\bB^{\sbOmg}_{\mathcal{S}}\big)^{1/2}
	\end{array}\right).
	\end{align*}
\end{proposition}

We see from Proposition \ref{prop1} that when $\bOmg = \bOmg_0$, it holds that $\bH^{\sbOmg_0} = \bI_{p+|\mathcal{S}|}$ and thus the value of the FDR function at point $\bOmg_0$ reduces to
\begin{align*}
{\rm FDR}_n(\bOmg_0, \mathcal{S})	= \mb{E}\Big[g_{n}\big(\bX_{\rm aug}^\mathcal{S}\big)\Big],
\end{align*}
which can be shown to be bounded from above by the target FDR level $q$ using the results proved in \cite{CandesFanJansonLv2016}. Since the dependence of FDR function on $\bOmg$ is completely through matrix $\bH^{\sbOmg}$, we can reparameterize the FDR function as ${\rm FDR}_n(\bH^{\sbOmg}, \mathcal{S})$.  In view of \eqref{eq: gn-def}, ${\rm FDR}_n(\bH^{\sbOmg}, \mathcal{S})$ is the expectation of some measurable function with respect to the probability law of $\bX_{\rm aug}^\mathcal{S}$ which has matrix normal distribution with independent rows, and thus is expected to be a smooth function of entries of $\bH^{\sbOmg}$ by measure theory. Motivated by such an observation, we make the following Lipschitz continuity assumption.

\begin{assumption}\label{cond0}
There exists some constant $L>0$ such that for all $|\mathcal{S}|\leq K_n$ and $\|\bOmg - \bOmg_{0}\|_2 \leq C_2a_n$ with some constant $C_2 > 0$ and $a_n\rightarrow 0$,
$\big|{\rm FDR}_n(\bH^{\sbOmg}, \mathcal{S}) - {\rm FDR}_n(\bH^{\sbOmg_0}, \mathcal{S})\big| \leq L\left\| \bH^{\sbOmg} - \bH^{\sbOmg_0}\right\|_F$,
	where $\|\cdot\|_2$ and $\|\cdot \|_F$ denote the matrix spectral norm and matrix Frobenius norm, respectively.
\end{assumption}

\begin{assumption}\label{cond0-1}
Assume that the estimated precision matrix $\wh{\bOmg}$  satisfies $\|\wh{\bOmg} - \bOmg_{0}\|_2 \leq C_2a_n$ with probability $1-O(p^{-c_1})$  for some constants $C_2, c_1 > 0$ and $a_n\rightarrow 0$, and that
$|\wt{\mathcal{S}}|\leq K_n$.
\end{assumption}

The error rate of precision matrix estimation assumed in Condition \ref{cond0-1} is quite flexible. We would like to emphasize that no sparsity assumption has been made on the true precision matrix $\bOmg_{0}$. Bounding the size of sparse models is also important for ensuring model identifiability and stability; see, for instance, \cite{FanLv2013, Lv2013} for more detailed discussions.

\begin{theorem}\label{thm1}
Assume that all the eigenvalues of $\bOmg_0$ are bounded away from $0$ and $\infty$ and the smallest eigenvalue of $2\diag(\bs) - \diag(\bs)\bOmg_0\diag(\bs)$ is bounded from below by some positive constant. Then under Condition \ref{cond0}, it holds that
	\begin{equation}\label{eq: e042}
	\sup_{|\mathcal{S}|\leq K_n,\, \|\sbOmg - \sbOmg_0\|_2\leq C_2 a_n}|{\rm FDR}_n(\bH^{\sbOmg}, \mathcal{S}) - {\rm FDR}(\bH^{\sbOmg_0}, \mathcal{S})| \leq O(K_n^{1/2} a_n).
	\end{equation}
Moreover under Conditions \ref{cond0}--\ref{cond0-1} with  $K_n^{1/2} a_n\rightarrow 0$, the FDR of \emph{RANK} is bounded from above by $q + O(K_n^{1/2} a_n) + O(p^{-c_1})$,
where $q \in (0, 1)$ is the target FDR level.
\end{theorem}

Theorem \ref{thm1} establishes the robustness of the FDR with respect to the precision matrix $\bOmg$; see the uniform bound in (\ref{eq: e042}). As a consequence, it shows that our modified model-free knockoffs procedure RANK can indeed have FDR asymptotically controlled at the target level $q$. We remark that the term $K_n^{1/2}$ in Theorem \ref{thm1} is because Condition \ref{cond0} is imposed through the matrix Frobenius norm, which is motivated from results on the smoothness of integral function from calculus. If one is willing to impose assumption through matrix spectral norm instead of Frobenius norm, then the extra term $K_n^{1/2}$ can be dropped and the set $\mathcal{S}$ can be taken as the full model $\{1,\cdots, p\}$.

We would like to stress that Theorem \ref{thm1} allows for arbitrarily complicated dependence structure of response $\by$ on covariates $\bx$ and for any valid construction of knockoff statistics $W_j$'s. This is different from the conditions needed for power analysis in Section \ref{sec: oracle-power} (that is, the linear model setting and LCD knockoff statistics). Moreover, the asymptotic FDR control in Theorem \ref{thm1} does not need the sure screening property of $\mathbb{P}\{\wt{\mathcal{S}} \supset \mathcal{S}_0\} \rightarrow 1$ as $n \rightarrow \infty$.

\subsection{Robustness of power in linear models}

We are now curious about the other side of the coin; that is, the robustness theory for the power of our modified knockoffs procedure RANK. As argued at the beginning of Section \ref{sec: oracle-power}, to ease the presentation and simplify the technical derivations we come back to high-dimensional linear models (\ref{eq001}) and use the LCD in (\ref{WH-statistics0}) as the knockoff statistics.  The difference with the setting in Section \ref{sec: oracle-power} is that we no longer assume that the true precision matrix $\bOmg_{0}$ is known and use the modified knockoffs procedure introduced in Section \ref{sec: mod-knockoffs} to achieve asymptotic FDR control.

Recall that for the RANK procedure,  the reduced model $\wt{\mathcal{S}}$ is first obtained from an independent subsample and then the knockoffs procedure is applied on the second fold of data to further select features from $\wt{\mathcal{S}}$. Clearly if $\wt{\mathcal{S}}$ does not have the sure screening property of $\mathbb{P}\{\wt{\mathcal{S}} \supset \mathcal{S}_0\} \rightarrow 1$ as $n \rightarrow \infty$, then the Lasso solution based on $[\bX_{\wt{\mathcal{S}}}^{(2)}, \tbX_{\wt{\mathcal{S}}}^{\sbOmg}]$ as given in \eqref{SLasso-solution}  is no longer a consistent estimate of $\bbeta_0$ even when the true precision matrix $\bOmg_0$ is used to generate the knockoff variables. In addition, the final power of our modified knockoffs procedure will always be upper bounded by $s^{-1}|\wt{\mathcal{S}} \cap \mathcal{S}_0|$. Nevertheless, the results in this section are still useful in the sense that model \eqref{eq001} can be viewed as the projected model on support $\wt{\mathcal{S}}$. Thus our power analysis here is relative power analysis with respect to the reduced model $\wt{\mathcal{S}}$. In other words, we will focus on how much power loss would occur after we apply the model-free knockoffs procedure to $(\bX_{\wt{\mathcal{S}}}^{(2)}, \tbX_{\wt{\mathcal{S}}}^{{\sbOmg}}, \by^{(2)})$ when compared to the power of $s^{-1}|\wt{\mathcal{S}} \cap \mathcal{S}_0|$. Since our focus is relative power loss, without loss of generality we will condition on the event
\begin{equation}\label{event-screening}
\left\{\wt{\mathcal{S}} \supset \mathcal{S}_0\right\}.
\end{equation}
We would like to point out that all conditions and results in this section can be adapted correspondingly when we view model \eqref{eq001} as the projected model if $\wt{\mathcal{S}}\not\supset \mathcal{S}_0$. Similarly as in FDR analysis, we restrict ourselves to sparse models with size bounded by $K_n$ that diverges as $n \rightarrow \infty$, that is, $|\wt{\mathcal{S}}| \leq K_n$.


With $\bOmg$ taken as the estimated precision matrix $\wh{\bOmg}$, we can generate the knockoff variables from \eqref{eq: xtilde-distr}. Then the Lasso procedure can be applied to the augmented data  $(\bX^{(2)}, \hbX, \by^{(2)})$ with $\hbX$ constructed in Step 3 of our modified knockoffs  procedure and  the LCD can be defined as
\begin{eqnarray}\label{WH-statistics}
\wh{W}_j= W_j^{\wh{\sbOmg}, \wt{\mathcal{S}}}= |\widehat{\beta}_j(\lambda; \wh{\bOmg}, \wt{\mathcal{S}})|-|\widehat{\beta}_{p+j}(\lambda;\wh{\bOmg}, \wt{\mathcal{S}})|,
\end{eqnarray}
where $\widehat{\beta}_j(\lambda; \wh{\bOmg},\wt{\mathcal{S}})$ and $\widehat{\beta}_{p+j}(\lambda; \wh{\bOmg},\wt{\mathcal{S}})$ are the $j$th and $(j+p)$th components, respectively, of the Lasso estimator
\begin{equation}\label{SLasso-solution}
\widehat{\bbeta}(\lambda; \wh{\bOmg},\wt{\mathcal{S}})  ={\argmin}_{\bb_{\wt{\mathcal{S}}_1} = \bzero}\left\{n^{-1} \big\|\by^{(2)}-[\bX^{(2)}, \wh{\bX}]\bb\big\|_2^2+\lambda\|\bb\|_1\right\}
\end{equation}
with  $\lambda\geq 0$  the regularization parameter and $\wt{\mathcal{S}}_1 = \{1\leq j\leq 2p: j\not\in \wt{\mathcal{S}} \text{ and }
j-p \not\in  \wt{\mathcal{S}}\} $.



Unlike the FDR analysis in Section \ref{sec: 2.2}, we now need sparsity assumption on the true precision matrix $\bOmg_{0}$.

\begin{assumption}\label{cond2}
Assume that $\bOmg_{0}$ is $L_p$-sparse with each row having at most $L_p$ nonzeros for some diverging $L_p$ and all the eigenvalues of $\bOmg_{0}$ are bounded away from 0 and $\infty$.
\end{assumption}

For each given precision matrix $\bOmg$ and reduced model $\mathcal{S}$, we define $W_j^{\sbOmg, \mathcal{S}}$ similarly as in \eqref{WH-statistics} except that $\bOmg$ is used to generate the knockoff variables and set $\mathcal{S}$ is used in \eqref{SLasso-solution} to calculate the Lasso solution. Denote by $\wh{\mathcal{S}}^{\sbOmg} = \{j: W_j^{\sbOmg, \mathcal{S}} \geq T \} \subset \mathcal{S}$ the final set of selected features using the LCD $W_j^{\sbOmg, \mathcal{S}}$ in the knockoffs inference framework. We further define a class of precision matrices $\bOmg \in \mathbb{R}^{p\times p}$
\begin{equation} \label{neweq007}
\calA = \left\{\bOmg: \ \bOmg \text{ is } L_p^{\prime}\text{-sparse and } \|\bOmg - \bOmg_0\|_2 \leq C_2a_n\right\},
\end{equation}
where $C_2$ and $a_n$ are the same as in Theorem \ref{thm1} and $L_p^{\prime}$ is some positive integer that diverges with $n$.  Similarly as in Section \ref{sec: oracle-power}, in the technical analysis we assume implicitly that with asymptotic probability one, for all valid constructions of the knockoff variables there are no ties in the magnitude of nonzero knockoff statistics 
and no ties in the magnitude of nonzero components of Lasso solution uniformly over all $\bOmg \in \calA$ and $|\mathcal{S}|\leq K_n$.

\begin{assumption}\label{cond5}
It holds that $\mathbb{P}\{\wh{\bOmg} \in \calA\} = 1 - O(p^{-c_2})$ for some constant $c_2 > 0$.
\end{assumption}

The assumption on the estimated precision matrix $\wh{\bOmg}$ made in Condition \ref{cond5} is mild and flexible. A similar class of precision matrices was considered in \cite{FanKongLiZheng2015} with detailed discussions on the choices of the estimation procedures. See, for example, \cite{FanLv2016, ChenRenZhaoZhou2016} for some more recent developments on large precision matrix estimation. In parallel to Theorem \ref{theo1}, we have the following results on the power of our modified knockoffs procedure with the estimated precision matrix $\wh{\bOmg}$.

\begin{theorem} \label{theo2}
Assume that Conditions \ref{cond4}--\ref{cond3} and \ref{cond2}--\ref{cond5} hold, the smallest eigenvalue of $2\diag(\bs) - \diag(\bs)\bOmg_0\diag(\bs)$ is positive and bounded away from $0$, $|\{j: |\beta_{0,j}|\gg [sn^{-1}(\log p)]^{1/2}\}| \geq cs$, and $\lambda =  C_3 \{(\log p)/n\}^{1/2}$ with $ c  \in ((qs)^{-1},1)$ and $C_3>0$ some constants. Then if $[(L_p + L_p^{\prime})^{1/2}+K_n^{1/2}] a_n = o(1)$ and $s \{a_n + (K_n+L_{p}') [n^{-1}(\log p)]^{1/2} \} = o(1)$, \emph{RANK} with estimated precision matrix $\wh{\bOmg}$ and reduced model $\wt{\mathcal S}$ has asymptotic power one.
\end{theorem}

Theorem \ref{theo2} establishes the robustness of the power for the RANK method. In view of Theorems \ref{thm1}--\ref{theo2}, we see that our modified knockoffs procedure RANK can enjoy appealing properties of FDR control and power simultaneously when the true covariate distribution is unknown and needs to be estimated in high dimensions.

\section{Simulation studies} \label{sec4}
So far we have seen that our suggested RANK method admits appealing theoretical properties for large-scale inference in high-dimensional nonlinear models. 
We now examine the finite-sample performance of RANK through four simulation examples.

\subsection{Model setups and simulation settings} \label{sec4.1}
Recall that the original knockoff filter (KF) in \cite{BarberCandes2015} was designed for linear regression model with dimensionality $p$ not exceeding sample size $n$, while the high-dimensional knockoff filter (HKF) in \cite{BarberCandes2016} considers linear model with $p$ possibly larger than $n$. To compare RANK with the HKF procedure in high-dimensional setting, 
our first simulation example adopts the linear regression model
\begin{equation} \label{simu-lm-1}
\by = \bX\bbeta + \bveps,
\end{equation}
where $\by$ is an $n$-dimensional response vector, $\bX$ is an $n \times p$ design matrix, $\bbeta = (\beta_1, \cdots, \beta_p)\t$ is a $p$-dimensional regression coefficient vector, and $\bveps$ is an $n$-dimensional error vector. Nonlinear models provide useful and flexible alternatives to linear models and are widely used in real applications. Our second through fourth simulation examples are devoted to three popular nonlinear model settings: the partially linear model, the single-index model, and the additive model, respectively. As a natural extension of linear model (\ref{simu-lm-1}), the partially linear model assumes that
\begin{equation} \label{simu-plm-2}
\by = \bX\bbeta + {\bg}({\bU}) + \bveps,
\end{equation}
where $ {\bg}({\bU})=(g(U_1),\cdots,g(U_n))^T$ is an $n$-dimensional vector-valued function with covariate vector ${\bU} = (U_1, \cdots, U_n)^T$, $g(\cdot)$ is some \textit{unknown} smooth nonparametric function, and the rest of notation is the same as in model (\ref{simu-lm-1}). In particular, the partially linear model is a semiparametric regression model that has been commonly used in many areas such as economics, finance, medicine, epidemiology, and environmental science \cite{Engle1986, Hardle2000}.

The third and fourth simulation examples drop the linear component. As a popular tool for dimension reduction, the single-index model assumes that
\begin{equation} \label{simu-sim-3}
\by = {\bg}(\bX\bbeta) + \bveps,
\end{equation}
where ${\bg}(\bX\bbeta)=(g(\bx_1^T\bbeta),\cdots,g(\bx_n^T\bbeta))^T$ with $\bX = (\bx_1, \cdots, \bx_n)\t$, $g(\cdot)$ is an \textit{unknown} link function, and the remaining notation is the same as in model (\ref{simu-lm-1}). In particular, the single-index model provides a flexible extension of the GLM by relaxing the parametric form of the link function \cite{ichimura1993semiparametric, stoker1986consistent, hardle1989investigating, li2007nonparametric, horowitz2009semiparametric}. To bring more flexibility while alleviating the curse of dimensionality, the additive model assumes that
\begin{equation} \label{simu-gam-4}
\by =\sum_{j=1}^p \bg_j(\bX_j) + \bveps,
\end{equation}
where $\bg_j(\btheta) = (g_j(\theta_1), \cdots, g_j(\theta_n))^T$ for $\btheta = (\theta_1, \cdots, \theta_n)\t$, $\bX_j$ represents the $j$th covariate vector with $\bX = (\bX_1, \cdots, \bX_p)$, $g_j(\cdot)$'s are some \textit{unknown} smooth functions, and the rest of notation is the same as in model (\ref{simu-lm-1}). The additive model has been widely employed for nonparametric modeling of high-dimensional data \cite{Hastie1990,
Ravikumar2009,Meier2009,
Chouldechova2015}.

For the linear model (\ref{simu-lm-1}) in simulation example 1, the rows of the $n \times p$ design matrix $\bX$ are generated as i.i.d. copies of $N(\bzero, \bSig)$ with covariance matrix $\bSig = (\rho^{|j - k|})_{1 \leq j, k \leq p}$ for $\rho = 0$ and $0.5$. We set the true regression coefficient vector $\bbeta_0\in \mathbb{R}^p$ as a sparse vector with $s = 30$ nonzero components, where the signal locations are chosen randomly and each nonzero coefficient is selected randomly from $\{\pm A\}$ with $A=1.5$ and $3.5$. The error vector $\bveps$ is assumed to be $N(\bzero, \sigma^2 \bI_n)$ with $\sigma = 1$. We set sample size $n = 400$ and consider the high-dimensional scenario with dimensionality $p = 200, 400, 600, 800$, and $1000$. For the partially linear model (\ref{simu-plm-2}) in simulation example 2, we choose the true function as $g(U)=\sin(2\pi U)$, generate ${\bU} = (U_1, \cdots, U_n)^T$ with i.i.d. $U_i$ from uniform distribution on $[0, 1]$, and set $A = 1.5$ with the remaining setting the same as in simulation example 1.

Since the single-index model and additive model are more complex than the linear model and partially linear model, we reduce the true model size $s$ while keeping sample size $n= 400$ in both simulation examples 3 and 4. For the single-index model (\ref{simu-sim-3}) in simulation example 3, we consider
the true link function $g(x) = x^3 / 2 $ and set $p = 200, 400, 600, 800$, and $1000$.
The true $p$-dimensional regression coefficient vector $\bbeta_0$ is generated similarly with $s = 10$ and $A = 1.5$. For the additive model (\ref{simu-gam-4}) in simulation example 4, we assume that $s = 10$ of the functions $g_j(\cdot)$'s are nonzero with $j$'s chosen randomly from $\{1, \cdots, p\}$ and the remaining $p - 10$ functions $g_j(\cdot)$'s vanish. Specifically, each nonzero function $g_j(\cdot)$ is taken to be a polynomial of degree 3 and all coefficients under the polynomial basis functions are generated independently 
as $N(0, 10^2)$ as in \cite{Chouldechova2015}.
The dimensionality $p$ is allowed to vary with values $200, 400, 600, 800$, and $1000$. For each simulation example, we set the number of repetitions as $100$.

\subsection{Estimation procedures} \label{sec4.2}

To implement RANK procedure described in Section \ref{sec: mod-knockoffs}, we need to construct a precision matrix estimator $\wh\bOmg$ and obtain the reduced model $\wt{\mathcal S}$ using the first fold of data $(\bX^{(1)}, \by^{(1)})$. Among all available estimators in the literature, we employ the ISEE method in \cite{FanLv2016} for precision matrix estimation due to its scalability, simple tuning, and nice theoretical properties. For simplicity, we choose $s_j=1/\Lambda_{\max}(\wh{\bOmg})$ for all $1 \leq j \leq p$, where $\wh{\bOmg}$ denotes the ISEE estimator for the true precision matrix $\bOmg_0$ and $\Lambda_{\max}$ standards for the largest eigenvalue of a matrix. Then we can obtain an $(n/2) \times (2p)$ augmented design matrix $[\bX^{(2)}, {\wh{\bX}}]$, where ${\wh{\bX}}$ represents an $(n/2) \times p$ knockoffs matrix constructed in Step 3 of our modified knockoffs procedure in Section \ref{sec: mod-knockoffs}. 
To construct the reduced model $\wt{\mathcal S}$ using the first fold of data $(\bX^{(1)}, \by^{(1)})$, we borrow the strengths from the recent literature on feature selection methods. After $\wt{\mathcal S}$ is obtained, we employ the reduced data $(\bX_{\rm aug}^{\wt{\mathcal S}}, \by^{(2)})$ with $\bX_{\rm aug}^{\wt{\mathcal S}} = [\bX_{\wt{\mathcal S}}^{(2)}, {\wh{\bX}}_{\wt{\mathcal S}}]$ to fit a model and construct the knockoff statistics.
In what follows, we will discuss feature selection methods for obtaining $\wt{\mathcal{S}}$ for the linear model (\ref{simu-lm-1}), partially linear model (\ref{simu-plm-2}), single-index model (\ref{simu-sim-3}), and additive model (\ref{simu-gam-4}) in simulation examples 1--4, respectively. We will also discuss the construction of knockoff statistics in each model setting. 

For the linear model (\ref{simu-lm-1}) in simulation example 1, we obtain the reduced model $\wt{\mathcal S}$ by first applying the Lasso procedure
\begin{equation} \label{neweq001}
\hbbeta^{(1)}= \argmin_{\bb\in \mathbb{R}^{p}} \left\{n^{-1} \|\by^{(1)} - \bX^{(1)} \bb\|_2^2 + \lambda \|\bb\|_1\right\}
\end{equation}
with $\lambda \geq 0$ the regularization parameter and then taking the support $\wt{\mathcal S} = \supp(\hbbeta^{(1)})$. 
Then with the estimated $\wh\bOmg$ and $\wt{\mathcal S}$, we construct the knockoff statistics as the LCD \eqref{WH-statistics}, where the estimated regression coefficient vector is obtained by applying the Lasso procedure on the reduced model as described in \eqref{SLasso-solution}. 
The regularization parameter $\lambda$ in Lasso is tuned using the $K$-fold cross-validation (CV).

For the partially linear model (\ref{simu-plm-2}) in simulation example 2, we employ the profiling method in semiparametric regression based on the first fold of data $(\bX^{(1)}, \bU^{(1)}, \by^{(1)})$ by observing that model (\ref{simu-plm-2}) becomes a linear model when conditioning on the covariate vector ${\bU}^{(1)}$. Consequently we need to estimate both the profiled response $\mathbb{E}(\by^{(1)} | \bU^{(1)})$ and the profiled covariates $\mathbb{E}(\bX^{(1)} | \bU^{(1)})$. To this end, we adopt the local linear smoothing estimators \cite{FanGijbels1996} $\widehat{\mathbb{E}(\by^{(1)} | \bU^{(1)})}$ and $\widehat{\mathbb{E}(\bX^{(1)} | \bU^{(1)})}$ of $\mathbb{E}(\by^{(1)} | \bU^{(1)})$ and $\mathbb{E}(\bX^{(1)} | \bU^{(1)})$ using the Epanechnikov kernel $K(u) = 0.75(1-u^2)_{+}$ with the optimal bandwidth selected by the generalized cross-validation (GCV). Then we define the Lasso estimator $\hbbeta^{(1)}$ for the $p$-dimensional regression coefficient vector similarly as in (\ref{neweq001}) with $\by^{(1)}$ and $\bX^{(1)}$ replaced by $\by^{(1)} - \widehat{\mathbb{E}(\by^{(1)} | \bU^{(1)})}$ and $\bX^{(1)} - \widehat{\mathbb{E}(\bX^{(1)} | \bU^{(1)})}$, respectively. The reduced model is then taken as $\wt{\mathcal S} = \supp(\hbbeta^{(1)})$. For knockoff statistics $\wh{W}_j$, we set $\wh{W}_j = 0$ for all $j \not\in \wt{\mathcal S}$. On the support $\wt{\mathcal S}$, we construct $\wh{W}_j = |\hbeta_{j}| - |\hbeta_{p+j}|$ with $\hbeta_{j}$ and $\hbeta_{p+j}$ the Lasso coefficients obtained by applying the model fitting procedure described above to the reduced data $(\bX_{\rm aug}^{\wt{\mathcal S}}, \bU_{\wt{\mathcal S}}^{(2)}, \by^{(2)})$ in the second subsample with  $ \bX_{\rm aug}^{\wt{\mathcal S}} = [\bX_{\wt{\mathcal S}}^{(2)}, {\wh{\bX}}_{\wt{\mathcal S}}]$.

To fit the single-index model (\ref{simu-sim-3}) in simulation example 3, we employ the Lasso-SIR method in \cite{LinZhaoLiu2016}. The Lasso-SIR first divides the sample of $m=n/2$ observations in the first subsample $(\bX^{(1)}, \by^{(1)})$  into $H$ slices of equal length $c$, and constructs $ \Lambda _{H} = \dfrac{1}{mc} (\bX^{(1)})^{T} \bM \bM ^{T} \bX^{(1)}$, where $ \bM = \bI_{H} \otimes \bone_{c}$ is an $m \times H$ matrix
that is the Kronecker product of the identity matrix $\bI_{H}$ and the constant vector $\bone_{c}$ of ones. Then the Lasso-SIR estimates the $p$-dimensional regression coefficient vector $\hbbeta^{(1)}$ using the Lasso procedure similarly as in (\ref{SLasso-solution}) with the original response vector $\by^{(1)}$ replaced by a new response vector $\widetilde{\by}^{(1)} = (c\lambda_1)^{-1} \bM \bM ^{T} \bX^{(1)} \bleta_1$, where $\lambda_1$ denotes the largest eigenvalue of matrix $\Lambda _{H}$ and $\bleta_1$ is the corresponding eigenvector. We set the number of slices $H = 5$. Then the reduced model is taken as $\wt{\mathcal S} = \supp(\hbbeta^{(1)})$. We then apply the fitting procedure Lasso-SIR discussed above  to the reduced data $(\bX_{\rm aug}^{\wt{\mathcal S}},  \by^{(2)})$ with $ \bX_{\rm aug}^{\wt{\mathcal S}} = [\bX_{\wt{\mathcal S}}^{(2)}, {\wh{\bX}}_{\wt{\mathcal S}}]$ and construct knockoff statistics in a similar way as in partially linear model.

 To fit the additive model (\ref{simu-gam-4}) in simulation example 4, we apply the GAMSEL procedure in \cite{Chouldechova2015} for sparse additive regression. In particular, we choose 6 basis functions each with 6 degrees of freedom for the smoothing splines using orthogonal polynomials
for each additive component and set the penalty mixing parameter $\gamma= 0.9$ in GAMSEL to obtain estimators of the true functions $g_j(\cdot)$'s. The GAMSEL procedure is first applied to the first subsample $(\bX^{(1)}, \by^{(1)})$ to obtain the reduced model $\wt{\mathcal S}$, and then applied to the reduced data   $(\bX_{\rm aug}^{\wt{\mathcal S}},  \by^{(2)})$ with $ \bX_{\rm aug}^{\wt{\mathcal S}} = [\bX_{\wt{\mathcal S}}^{(2)}, {\wh{\bX}}_{\wt{\mathcal S}}]$ to obtain estimates  $\widehat{g}_j$ and $\widehat{g}_{p+j}$ for the additive functions corresponding to the $j$th covariate and its knockoff counterpart with $j\in \wt{\mathcal {S}}$,  respectively. The  knockoff statistics are then constructed as
\begin{equation} \label{neweq003}
\wh{W}_j = \|\widehat{g}_j\|_{n/2}^2-\|\widehat{g}_{p+j}\|_{n/2}^2 \text{ for } j\in \wt{\mathcal S}
\end{equation}
and $\wh{W}_j  = 0$ for $j\not\in \wt{\mathcal S}$, where $\|\widehat{g}_j\|_{n/2}$ represents the empirical norm of the estimated function $\widehat{g}_j(\cdot)$ evaluated at its observed points and $n/2$ stands for the size of the second subsample.
%

It is seen that in all four examples above, intuitively large positive values of knockoff statistics $\wh W_j$ provide strong evidence against the $j$th null hypothesis $H_{0, j}: \beta_j = 0$ or $H_{0, j}: g_j = 0$. For all simulation examples, we set the target FDR level at $q = 0.2$.

\begin{table}[htp!]
	\centering
	\caption{\label{tab1} Simulation results for linear model (\ref{simu-lm-1}) in simulation example 1 with $A=1.5$ in Section \ref{sec4.1}}
	\vskip0.2cm
	\tabcolsep 0.15cm
	\begin{tabular}{ccccccccccccc}
		\toprule
	&	&  \multicolumn{2}{c}{RANK}&&\multicolumn{2}{c}{RANK$_+$} &&\multicolumn{2}{c}{RANKs}&&\multicolumn{2}{c}{RANKs$_+$}     \\
		\cline{3-4}\cline{6-7}\cline{9-10}\cline{12-13}
$\rho$	&	$ p$ &FDR & Power &  & FDR& Power && FDR & Power &  & FDR& Power  \\
		\hline
0&		200 & 0.2054 & 1.00 && 0.1749 & 1.00 &&0.1909  & 1.00   && 0.1730 & 1.00 \\
	&	400 & 0.2062 & 1.00 && 0.1824 & 1.00 && 0.2010 & 1.00   && 0.1801 & 1.00 \\
	&	600 & 0.2263 & 1.00 &&0.1940  & 1.00 && 0.2206 & 1.00   && 0.1935 & 1.00 \\
	&	800 & 0.2385 & 1.00 &&0.1911  & 1.00 &&0.2247  & 1.00   &&0.1874  & 1.00 \\
        & 1000    & 0.2413 & 1.00 &&0.2083  & 1.00 &&0.2235  & 1.00   &&0.1970  & 1.00 \\
0.5   &    200   & 0.2087 & 1.00 &&0.1844  & 1.00 && 0.1875 & 1.00   && 0.1692 & 1.00  \\
	&	400 & 0.2144 & 1.00 &&0.1879  & 1.00 &&0.1954  &1.00    && 0.1703 & 1.00 \\
	&	600 & 0.2292 & 1.00 &&0.1868  & 1.00 &&0.2062  & 1.00   && 0.1798 & 1.00  \\
	&	800 & 0.2398 & 1.00 &&0.1933  & 1.00 &&0.2052  &0.9997&& 0.1805 &0.9997 \\
        & 1000    & 0.2412 & 1.00 &&0.2019  & 1.00 &&0.2221  &0.9984&& 0.2034 &0.9984 \\
		\bottomrule
	\end{tabular}
\end{table}

\begin{table}[htp!]
		\centering
		\caption{\label{tab3} Simulation results for linear model (\ref{simu-lm-1}) in simulation example 1 with $A=3.5$ in Section \ref{sec4.1}}
		\vskip0.2cm
		\tabcolsep 0.15cm

	\begin{tabular}{ccccccccccccc}
		\toprule
&		&  \multicolumn{2}{c}{RANKs}&&\multicolumn{2}{c}{RANKs$_+$} &&\multicolumn{2}{c}{HKF}&&\multicolumn{2}{c}{HKF$_+$}     \\
		\cline{3-4}\cline{6-7}\cline{9-10}\cline{12-13}
$\rho$	&	$ p$ &FDR & Power &  & FDR& Power && FDR & Power &  & FDR& Power  \\
		\hline
0 &			200  &0.1858 &1.00     &&0.1785 &1.00     &&0.1977  &0.9849  &&0.1749  &0.9837 \\
		&	400  &0.1895 &1.00     && 0.1815&1.00     &&0.2064  &0.9046  &&0.1876  &0.8477 \\
		&	600  &0.2050 &1.00     &&0.1702 &1.00     &&0.1964  &0.8424  &&0.1593  &0.7668 \\
		&	800  &0.2149 &1.00     &&0.1921 &1.00     &&0.1703  &0.7513  &&0.1218  &0.6241 \\
                & 1000      &0.2180 &1.00     &&0.1934 &1.00     &&0.1422  &0.7138  &&0.1010  &0.5550 \\
0.5&        200            &0.1986 &1.00     && 0.1618&1.00     &&0.1992  &0.9336  &&0.1801  &0.9300 \\
		&	400  &0.1971 &1.00     &&0.1805 &1.00     &&0.1657  & 0.8398 &&0.1363  & 0.7825 \\
		&	600  &0.2021 &1.00     &&0.1757 &1.00     &&0.1253  & 0.7098 &&0.0910  &0.6068  \\
		&	800  &0.2018 &1.00     &&0.1860 &1.00     &&0.1374  &0.6978  &&0.0917  &0.5792 \\
                   & 1000   &0.2097 &0.9993 &&0.1920 &0.9993 &&0.1552  &0.6486  &&0.1076  &0.5524 \\
		\bottomrule
	\end{tabular}
	\end{table}

\begin{table}[htp!]
	\centering
	\caption{\label{tab-plm} Simulation results for partially linear model (\ref{simu-plm-2}) in simulation example 2 in Section \ref{sec4.1}}
	\vskip0.2cm
	\tabcolsep 0.15cm
	\begin{tabular}{ccccccccccccc}
		\toprule
	&	&  \multicolumn{2}{c}{RANK}&&\multicolumn{2}{c}{RANK$_+$} &&\multicolumn{2}{c}{RANKs}&&\multicolumn{2}{c}{RANKs$_+$}     \\
		\cline{3-4}\cline{6-7}\cline{9-10}\cline{12-13}
$\rho$	&	$ p$ &FDR & Power &  & FDR& Power && FDR & Power &  & FDR& Power  \\
		\hline
0&		200 &0.2117 & 1.00  &&0.1923  &1.00  &&0.1846  &0.9976  &&0.1699  &0.9970           \\
	&	400 &0.2234 & 1.00  &&0.1977  &1.00  &&0.1944  &0.9970  &&0.1747  &0.9966           \\
	&	600 &0.2041 & 1.00  &&0.1776  & 1.00 &&0.2014  &0.9968  &&0.1802  &0.9960           \\
	&	800 &0.2298 & 1.00  &&0.1810  & 1.00 &&0.2085  &0.9933  &&0.1902  &0.9930           \\
        & 1000    &0.2322 &  1.00  &&0.1979 & 1.00 && 0.2113  &0.9860  &&0.1851  &0.9840           \\
0.5   &	200 &0.2180 &1.00   &&0.1929  &1.00  &&0.1825  &0.9952  &&0.1660  &0.9949       \\
	&	400 &0.2254 &1.00   &&0.1966  &1.00  &&0.1809  &0.9950  &&0.1628  &0.9948       \\
	&	600 &0.2062 &1.00   &&0.1814  &1.00  &&0.2038  &0.9945  &&0.1898  &0.9945       \\
	&	800 &0.2264 &1.00   &&0.1948  &1.00  &&0.2019  &0.9916  &&0.1703  &0.9906        \\
           & 1000 &0.2316 & 1.00  &&0.2033  &1.00  &&0.2127  &0.9830  &&0.1857  &0.9790           \\
		\bottomrule
	\end{tabular}
\end{table}

\begin{table}[htp!]
	\centering
	\caption{\label{tab-sim} Simulation results for single-index model (\ref{simu-sim-3}) in simulation example 3 in Section \ref{sec4.1}
	}
	\vskip0.2cm
	\tabcolsep 0.15cm
	\begin{tabular}{ccccccccccccc}
		\toprule
	&	&  \multicolumn{2}{c}{RANK}&&\multicolumn{2}{c}{RANK$_+$} &&\multicolumn{2}{c}{RANKs}&&\multicolumn{2}{c}{RANKs$_+$}     \\
		\cline{3-4}\cline{6-7}\cline{9-10}\cline{12-13}
$\rho$	&	$ p$ &FDR & Power &  & FDR& Power && FDR & Power &  & FDR& Power  \\
		\hline
0   & 200  & 0.1893 & 1 &  & 0.1413 & 1 &  & 0.1899 & 1     &  & 0.1383 & 1     \\
    & 400  & 0.2163 & 1 &  & 0.1598 & 1 &  & 0.245  & 0.998 &  & 0.1676 & 0.997 \\
    & 600  & 0.2166 & 1 &  & 0.1358 & 1 &  & 0.2314 & 0.999 &  & 0.1673 & 0.998 \\
    & 800  & 0.1964 & 1 &  & 0.1406 & 1 &  & 0.2443 & 0.992 &  & 0.1817 & 0.992 \\
    & 1000 & 0.2051 & 1 &  & 0.134  & 1 &  & 0.2431 & 0.969 &  & 0.1611 & 0.962 \\
0.5 & 200  & 0.2189 & 1 &  & 0.1591 & 1 &  & 0.2322 & 1     &  & 0.1626 & 1     \\
    & 400  & 0.2005 & 1 &  & 0.1314 & 1 &  & 0.2099 & 0.996 &  & 0.1615 & 0.995 \\
    & 600  & 0.2064 & 1 &  & 0.1426 & 1 &  & 0.2331 & 0.998 &  & 0.1726 & 0.998 \\
    & 800  & 0.2049 & 1 &  & 0.1518 & 1 &  & 0.2288 & 0.994 &  & 0.1701 & 0.994 \\
    & 1000 & 0.2259 & 1 &  & 0.1423 & 1 &  & 0.2392 & 0.985 &  & 0.185  & 0.983 \\
		\bottomrule
	\end{tabular}
\end{table}

\begin{table}[htp!]
	\centering
	\caption{\label{tab-gam} Simulation results for additive model (\ref{simu-gam-4}) in simulation example 4 in Section \ref{sec4.1} 
}
	\vskip0.2cm
	\tabcolsep 0.15cm
	\begin{tabular}{ccccccccccccc}
		\toprule
	&	&  \multicolumn{2}{c}{RANK}&&\multicolumn{2}{c}{RANK$_+$} &&\multicolumn{2}{c}{RANKs}&&\multicolumn{2}{c}{RANKs$_+$}     \\
		\cline{3-4}\cline{6-7}\cline{9-10}\cline{12-13}
$\rho$	&	$ p$ &FDR & Power &  & FDR& Power && FDR & Power &  & FDR& Power  \\
		\hline
0&		200 &0.1926 &0.9780 &&0.1719 &0.9690 &&0.2207 &0.9490 &&0.1668 & 0.9410 \\
	&	400 &0.2094 &0.9750 &&0.1773 &0.9670 &&0.2236 &0.9430 &&0.1639 & 0.9340 \\
	&	600 &0.2155 &0.9670 &&0.1729 &0.9500 &&0.2051 &0.9310 &&0.1620 & 0.9220 \\
	&	800 &0.2273 &0.9590 &&0.1825 &0.9410 &&0.2341 &0.9280 &&0.1905  &0.9200       \\
        &  1000   &0.2390 & 0.9570 &&0.1751 &0.9350 &&0.2350 &0.9140 &&0.1833 &0.9070           \\
0.5   &      200 &0.1904 &0.9680  &&0.1733 & 0.9590&&0.2078 &0.9370  &&0.1531  &0.9330    \\
      &	400         &0.2173 & 0.9650 &&0.1701 & 0.9540&&0.2224 &0.9360  &&0.1591  &0.9280           \\
&	        600 &0.2267 & 0.9600 &&0.1656 & 0.9360&&0.2366 &0.9340  &&0.1981  &0.9270           \\
&	        800 &0.2306 & 0.9540 &&0.1798 & 0.9320&&0.2332 &0.9150  &&0.1740  &0.9110           \\
&          1000   &0.2378 &0.9330 &&0.1793  &0.9270 &&0.2422 &0.8970  &&0.1813  &0.8880           \\
		\bottomrule
	\end{tabular}
\end{table}

\subsection{Simulation results} \label{sec4.3}
To gain some insights into the effect of data splitting, we also implemented our procedure without the data splitting step.
To differentiate, we use RANKs to denote the procedure with data splitting and RANK to denote the procedure without data splitting. To examine the feature selection performance, we look at both measures of FDR and power. The empirical versions of FDR and power based on 100 replications are reported in Tables \ref{tab1}--\ref{tab3} for simulation example 1 and Tables \ref{tab-plm}--\ref{tab-gam} for simulation examples 2--4, respectively. In particular, Table \ref{tab1} compares the performance of RANK and RANK$_+$ with that of RANKs and RANKs$_+$, where the subscript $+$ stands for the corresponding method when the modified knockoff threshold $T_+$ is used.   
We see from Table \ref{tab1} that RANK and RANK$_+$ mimic closely RANKs and RANKs$_+$, respectively, suggesting that data splitting is more of a technical assumption.  In addition, the FDR is approximately controlled at the target level of $q = 0.2$ with high power, which is in line with our theory. Table \ref{tab3} summarizes the comparison of RANKs with HKF procedure for high-dimensional linear regression model. Despite that both methods are based on data splitting, their practical performance is very different.  It is seen that although controlling the FDR below the target level, HKF suffers from a loss of power due to the use of the screening step and the power deteriorates as dimensionality $p$ increases. In contrast, the performance of RANKs is robust across different correlation levels $\rho$ and dimensionality $p$. It is worth mentioning that HKF procedure with data recycling performed generally better than that with data splitting alone. Thus only the results for the former version are reported in Table \ref{tab3} for simplicity.

For high-dimensional nonlinear settings of partially linear model, single-index model, and additive model in simulation examples 2--4, we see from Tables \ref{tab-plm}--\ref{tab-gam} that RANKs and RANKs$_+$ performed well and similarly as RANK and RANK$_+$  in terms of both FDR control and power across different scenarios.   These results demonstrate the  model-free feature of our procedure for large-scale inference in nonlinear models.  

\section{Real data analysis} \label{sec5}

In addition to simulation examples presented in Section \ref{sec4}, we also demonstrate the practical utility of our RANK procedure on a gene expression data set, which is based on Affymetrix GeneChip microarrays for the plant Arabidopsis thaliana  in \cite{Wille2004}. It is well known that isoprenoids play a key role in plant and animal physiological processes, such as  photosynthesis,  respiration, regulation of growth, and defense against pathogens in plant physiological processes.
In particular, \cite{Horvath2003} found that many of the genes expressed preferentially
in mature leaves  are readily recognizable as genes
involved in photosynthesis, including  rubisco activase (AT2G39730), fructose bisphosphate aldolase
(AT4G38970), and two glycine hydroxymethyltransferase
genes (AT4G37930 and AT5G26780).  Thus isoprenoids have become important ingredients in various drugs (e.g., against cancer and malaria), fragrances (e.g., menthol), and food colorants (e.g., carotenoids). See, for instance, \cite{Wille2004, Schafer2005, Prelic2006} on studying the mechasnism of isoprenoid synthesis in a wide range of applications.

The aforementioned data set in \cite{Wille2004} consists of $118$ gene expression patterns under various experimental conditions for $39$ isoprenoid genes, 15 of which are assigned to the regulatory pathway, 19 to the plastidal pathway, and the remaining 5 isoprenoid genes encode protein located in the mitochondrion.  Moreover, 795 additional genes from 56 metabolic pathways are incorporated into the isoprenoid genetic network.  Thus the combined data set is comprised of a sample of $n = 118$ gene expression patterns for $834$ genes.  This data set was studied in  \cite{Yang2014} for identifying genes that exhibit significant association with the specific isoprenoid gene GGPPS11 (AGI code AT4G36810).
Motivated by   \cite{Yang2014}, we choose the expression level of isoprenoid gene GGPPS11 as the response and treat the remaining $p = 833$ genes from 58 different metabolic pathways as the covariates, in which the dimensionality $p$ is much larger than sample size $n$. All the variables are logarithmically transformed. To identify important genes associated with isoprenoid gene GGPPS11, we employ the RANK method using the Lasso procedure with target FDR level $q=0.2$. The implementation of RANK is the same as that  in Section \ref{sec4} for the linear model. Since the sample size of this data set is relatively low, we choose to implement RANK without sample splitting, which has been demonstrated in Section \ref{sec4} to be capable of controlling the FDR at the desired level.



\begin{table}[tp!]
	\centering
	\caption{\label{tab5-real} Selected genes and their associated pathways for real data analysis in Section \ref{sec5}}
	\vskip0.2cm
	\tabcolsep 0.6cm
	\begin{tabular}{lcclc}
		\toprule
		\multicolumn{2}{c}{RANK}&&\multicolumn{2}{c}{RANK$_+$}      \\
		\cline{1-2}\cline{4-5}
		Pathway &Gene && Pathway &Gene    \\
		\hline
		Calvin   & AT4G38970  &   &   Calvin   & AT4G38970  \\
		Carote  & AT1G57770 &   &   Carote   & AT1G57770  \\
		Folate   & AT1G78670 &   &   Folate   & AT1G78670  \\
		Inosit     & AT3G56960 &   &      &  \\
		Phenyl   & AT2G27820 &   &   Phenyl   & AT2G27820  \\
		Purine    & AT3G01820 &   &   Purine    & AT3G01820  \\
		Ribofl     & AT4G13700 &   &   &  \\
		Ribofl     & AT2G01880 &   &   Ribofl     & AT2G01880  \\
		Starch    & AT5G19220  &   &   Starch    & AT5G19220  \\
		\hline
		\multicolumn{5}{c}{Lasso}\\
		\hline
		Pathway &Gene && Pathway &Gene \\
		\cline{1-2}\cline{4-5}
		Berber & AT2G34810  &  & Porphy  & AT4G18480   \\
		Calvin  & AT4G38970  &  &  Pyrimi  &  AT5G59440  \\
		Calvin  & AT3G04790  &  &  Ribofl   &  AT2G01880  \\
		Glutam & AT5G18170  &  &  Starch &  AT5G19220   \\
		Glycol   & AT4G27600  &  &  Starch &   AT2G21590 \\
		Pentos  & AT3G04790  &  &  Trypto &   AT5G48220 \\
		Phenyl   & AT2G27820  &  & Trypto  &   AT5G17980 \\
		Porphy  & AT1G03475  &  &  Mevalo &   AT5G47720 \\
		Porphy  & AT3G51820  &  &   &    \\
		\bottomrule
	\end{tabular}
\end{table}

Table \ref{tab5-real} lists the selected genes by RANK, RANK$_+$, and Lasso along with their associated pathways. 
We see from Table \ref{tab5-real} that RANK, RANK$_+$, and Lasso selected 9 genes, 7 genes, and 17 genes, respectively. The common set of four genes, AT4G38970, AT2G27820, AT2G01880, and AT5G19220, was selected by all three  methods. The values of the adjusted $R^2$ for these three selected models are equal to $0.7523$, $0.7515$, and $0.7843$, respectively, showing similar level of goodness of fit.
In particular, among the top 20 genes selected using the Elem-OLS method with entrywise transformed Gram matrix in \cite{Yang2014}, we found that five genes (AT1G57770, AT1G78670, AT3G56960, AT2G27820, and AT4G13700) selected by RANK are included in such a list of top 20 genes, and three genes (AT1G57770, AT1G78670,  and  AT2G27820) picked by RANK$_+$ are contained in the same list.

To gain some scientific insights into the selected genes, we conducted Gene Ontology (GO) enrichment analysis to interpret, from the biological point of view, the influence of selected genes on isoprenoid gene GGPPS11, 
which is known as a precursor to chloroplast, carotenoids, tocopherols, and abscisic acids. Specifically, in the enrichment test of GO biological process, gene AT1G57770 is involved in carotenoid biosynthetic process. In the GO cellular component enrichment test, 
genes AT4G38970 and AT5G19220 are located in chloroplast, chloroplast envelope, and chloroplast stroma; gene AT1G57770 is located in chloroplast and mitochondrion; and gene AT2G27820 is located in chloroplast, chloroplast stroma, and cytosol. The GO molecular function enrichment test shows that gene AT4G38970 has fructose-bisphosphate aldolase activity and gene AT1G57770 has carotenoid isomerase activity and oxidoreductase activity. These scientific insights in terms of biological process, cellular component, and molecular function suggest that the selected genes may have meaningful biological relationship with the target isoprenoid gene GGPPS11. See, for example, \cite{Horvath2003, Ramel2009, Wienkoop2004} for more discussions on these genes.


\section{Discussions} \label{sec6}

Our analysis in this paper reveals that the suggested RANK method exploiting the general framework of model-free knockoffs introduced in \cite{CandesFanJansonLv2016} can asymptotically control the FDR in general high-dimensional nonlinear models with \textit{unknown} covariate distribution. The robustness of the FDR control under \textit{estimated} covariate distribution is enabled by imposing the 
Gaussian graphical structure on the covariates. Such a structural assumption has been widely employed to model the association networks among the covariates and extensively studied in the literature. Our method and theoretical results are powered by scalable large precision matrix estimation with statistical efficiency. It would be interesting to extend the robustness theory of the FDR control beyond Gaussian designs as well as for heavy-tailed data and dependent observations.

Our work also provides a first attempt to the power analysis for the model-free knockoffs framework. The nontrivial technical analysis establishes that RANK can have asymptotic power one in high-dimensional linear model setting when the Lasso is used for sparse regression. It would be interesting to extend the power analysis for RANK with a wide class of sparse regression and feature screening methods including SCAD, SIS, and many other concave regularization methods \cite{FanLi2001, FanLv2008, FanFan2008, FanLv2013}. Though more challenging, it is also important to investigate the power property for RANK beyond linear models. Our RANK procedure utilizes the idea of data splitting, which plays an important role in our technical analysis. Our numerical examples, however, suggest that data splitting is more of a technical assumption than a practical necessity. It would be interesting to develop theoretical guarantees for  RANK without data splitting.  These extensions are interesting topics for future research.

\appendix
\section{Proofs of main results} \label{appA}

We provide the proofs of Theorems \ref{theo1}--\ref{theo2}, Propositions \ref{prop1}--\ref{pro1}, and Lemmas \ref{lem: cond-Shat0}--\ref{lemma: 1} in this appendix. Additional technical details for the proofs of Lemmas \ref{error-lem}--\ref{lem:3-prime} are included in the Supplementary Material. To ease the technical presentation, we first introduce some notation.
Let $\Lambda_{\min}(\cdot)$ and $\Lambda_{\max}(\cdot)$ be the smallest and largest eigenvalues of a symmetric matrix. For any matrix $\bA = (a_{ij})$, denote by $\|\bA\|_1 = \max_{j}\sum_{i} |a_{ij}|$, $\|\bA\|_{\max}= \max_{i,j}|a_{ij}|$, $\|\bA\|_2=\Lambda_{\max}^{1/2}(\bA^T \bA)$, and $\|\bA\|_F = [\tr(\bA^T \bA)]^{1/2}$ the matrix $\ell_1$-norm, entrywise maximum norm, spectral norm, and Frobenius norm, respectively. For any set $\mathcal{S} \subset \{1,\cdots, p\}$, we use $\bA_{\mathcal S}$ to represent the submatrix of $\bA$ formed by columns in set $\mathcal{S}$ and $\bA_{\mathcal{S},\mathcal{S}}$ to denote the principal submatrix formed by columns and rows in set $\mathcal{S}$.

\subsection{Proofs of Lemma \ref{lem: cond-Shat0} and Theorem \ref{theo1}} \label{appA.1}

Observe that the choice of $\wt{\mathcal{S}} = \{1, \cdots, p\}$ certainly satisfies the sure screening property. We see that Lemma \ref{lem: cond-Shat0} and Theorem \ref{theo1} are specific cases of Lemma \ref{lem: cond-Shat} in Section \ref{appB.4} of Supplementary Material and Theorem \ref{theo2}, respectively. Thus we only prove the latter ones.

\subsection{Proof of Proposition \ref{prop1}} \label{appA.2}

In this proof, we will consider $\bOmg$ and $\mathcal S$ as deterministic parameters and focus only on the second half of sample $(\bX^{(2)}, \by^{(2)})$ used in FDR control. Thus, we will drop the superscripts in $(\bX^{(2)}, \by^{(2)})$ whenever there is no confusion.
For a given precision matrix $\bOmg$, the matrix of knockoff variables
$$\wt{\bX}^{\sbOmg}= [\tbx_{1}^{\sbOmg}, \cdots, \tbx_{n}^{\sbOmg}]^T$$
can be generated using \eqref{eq: xtilde-distr} with $\bOmg_{0}$ replaced by  $\bOmg$. Here, we use the superscript $\bOmg$ to emphasize the dependence of knockoffs matrix on $\bOmg$.  Recall that for a given set $\mathcal{S}$ with $k=|\mathcal{S}|$, we calculate the knockoff statistics $W_j$'s using $(\by,\bX_S, \wt{\bX}_S^{\sbOmg})$. Thus, the FDR function can be written as
\begin{align}\label{eq: e005}
 {\rm FDR}_n(\bOmg, \mathcal{S})  = {\mathbb{E}}[{\rm FDP}_n(\by, \bX_\mathcal{S}, \wt{\bX}^{\sbOmg}_\mathcal{S})]
 = \mb{E}\big[g_{1,n}(\bX, \tbX_{\mathcal{S}}^{\sbOmg}) \big],
\end{align}
where
$
g_{1,n}(\bX, \tbX_{\mathcal{S}}^{\sbOmg})=\mathbb{E}[{\rm FDP}_n(\by, \bX_{\mathcal S}, \tbX_{\mathcal{S}}^{\sbOmg} )\big| \bX, \tbX_{\mathcal{S}}^{\sbOmg}]$.
It is seen that the function $g_n$ is the conditional FDP when knockoff variables $\tbx_i^{\sbOmg}$, $1 \leq i \leq n$, are simulated using $\bOmg$ and only variables in set $\mathcal{S}$ are used to construct knockoff statistics $W_j$. We want to emphasize that since given $\bX$ the response $\by$ is independent of $\tbX_{\mathcal{S}}^{\sbOmg}$, the functional form of $g_{1,n}$ is free  of the matrix $\bOmg$ used to generate knockoff variables.

Using the technical arguments in \cite{CandesFanJansonLv2016}, we can show that ${\rm FDR}_n(\bOmg_0, \mathcal{S}) \leq q$ for any sample size $n$ and all subsets $\mathcal{S} \subset \{1,\cdots, p\}$ that are independent of the original data $(\bX^{(2)}, \by^{(2)})$ used in the knockoffs procedure. Observe that the only difference between ${\rm FDR}_n(\bOmg, \mathcal{S})$ and ${\rm FDR}_n(\bOmg_0, \mathcal{S})$ is that different precision matrices are used to generate knockoff variables. We restrict ourselves to the following data generating scheme
\begin{align*}
 \tbx_i^{\sbOmg}  = (\bC^{\sbOmg})^T\bx_i+ \bB^{\sbOmg}\bz_i, \quad i=1,\cdots, n,
\end{align*}
where $\bC^{\sbOmg}=\bI_p- \bOmg\diag\{\bs\}$, $\bz_i \sim N(\bzero, \bI_{p})\in \mathbb{R}^{p}$ are i.i.d. normal random vectors that are independent of $\bx_i$'s, and  $ \bB^{\sbOmg}= \big(2\diag\{\bs\} - \diag\{\bs\}\bOmg\diag\{\bs\}\big)^{1/2}$. For simplicity, write $\tbx_i^{(0)} =  \tbx_i^{\sbOmg_0}$, $\bB_0=\bB^{\sbOmg_0}$, and $\bC_0=\bC^{\sbOmg_0}$, i.e., the matrices corresponding to the oracle case. Then restricted to set $\mathcal{S}$,
\begin{align*}
 \tbx^{\sbOmg}_{i,\mathcal{S}}  = (\bC_{\mathcal{S}}^{\sbOmg})^T\bx_i+ (\bB^{\sbOmg}_S)^T\bz_i,\quad \tbx^{(0)}_{i,\mathcal{S}}  = \bC_{0,\mathcal{S}}^T\bx_i+ (\bB_{0,\mathcal{S}})^T\bz_i,
\end{align*}
where the subscript $\mathcal{S}$ means the submatrix (subvector) formed by columns (components) in set $\mathcal{S}$.
We want to make connections between $ \tbx^{\sbOmg}_{i,\mathcal{S}}$ and $\tbx^{(0)}_{i,\mathcal{S}}$. To this end, 
construct
\begin{align}\label{eq: proxy-knockoff}
& \breve{\bx}^{\sbOmg}_{i,\mathcal{S}}  =(\bC_{\mathcal{S}}^{\sbOmg})^T\bx_i+ \wt{\bB}^T\bB_{0,\mathcal{S}}^T\bz_i,
\end{align}
where 
$\wt{\bB} = (\bB_{0,\mathcal{S}}^T\bB_{0,\mathcal{S}})^{-1/2}\big((\bB_{\mathcal{S}}^{\sbOmg})^T\bB_{\mathcal{S}}^{\sbOmg}\big)^{1/2}$.
Then it is seen that $(\bx_i, \tbx^{\sbOmg}_{i,\mathcal{S}})$ and $(\bx_i, \breve{\bx}^{\sbOmg}_{i,\mathcal{S}})$ have identical joint distribution. Although $\breve{\bx}^{\sbOmg}_{i,\mathcal{S}}$ cannot be calculated in practice for a given $\bOmg$ due to its dependency on $\bOmg_0$, the  random vector $\breve{\bx}^{\sbOmg}_{i,\mathcal{S}}$ acts as a proxy of $\tbx^{\sbOmg}_{i,\mathcal{S}}$ in studying the FDR function. In fact, by construction \eqref{eq: e005} can be further written as
\begin{align}\label{eq: FDR1}
 {\rm FDR}_n(\bOmg, \mathcal{S})
= \mb{E}\big[g_{1,n}(\bX, \tbX_{\mathcal{S}}^{\sbOmg}) \big]  = \mb{E}\big[g_{1,n}(\bX, \breve{\bX}_{\mathcal{S}}^{\sbOmg}) \big],
\end{align}
where $ \breve{\bX}_{\mathcal{S}}^{\sbOmg} = [\breve{\bx}^{\sbOmg}_{1,\mathcal{S}}, \cdots, \breve{\bx}^{\sbOmg}_{n,\mathcal{S}}]^T$.

Observe that the randomness in both $\tbX_{\mathcal{S}}^{(0)}$ and $\breve{\bX}_{\mathcal{S}}^{\sbOmg}$ is fully determined by the same random matrices $ \bX$ and $\bZ\bB_{0,\mathcal{S}}$, which are independent of each other and whose rows are i.i.d. copies from $N(\bzero,  \bSig_0)$ and $N(\bzero, \bB_{0,\mathcal{S}}^T\bB_{0,\mathcal{S}})$, respectively. 
For this reason, we can rewrite the FDR function in \eqref{eq: FDR1} as
\begin{align*}
& {\rm FDR}_n(\bOmg, \mathcal{S}) = \mb{E}[g_{1,n}( \bX_{\rm aug}^\mathcal{S}\bH^{\sbOmg})], 
\end{align*}
where $ \bX_{\rm aug}^\mathcal{S}  = [\bX, \tbX_{0,\mathcal{S}}]= [\bX, \bX\bC_{0,\mathcal S} + \bZ\bB_{0,\mathcal{S}}]  \in \mb{R}^{n\times (p+k)}$ is the augmented matrix collecting columns of $\bX$ and $\tbX_{0,\mathcal{S}}$, and
\begin{align*}
\bH^{\sbOmg} = \left(\begin{array}{cc}
\bI_{p} & \bC_{\mathcal S}^{\sbOmg}- \bC_{0,\mathcal{S}}\wt{\bB} \\
\bzero & \wt{\bB}
\end{array}\right),
\end{align*}
which completes the proof of Proposition \ref{prop1}.

\subsection{Lemma \ref{lemma: 1} and its proof} \label{appA.3}

\begin{lemma} \label{lemma: 1}
Assume that $\|\bOmg - \bOmg_0\|_2 = O(a_n)$ with $a_{n} \rightarrow 0$ some deterministic sequence and all the notation the same as in Proposition \ref{prop1}. If $\Lambda_{\min}\{2\diag(\bs) - \diag(\bs)\bOmg_0\diag(\bs)\} \geq c_0$ and $\Lambda_{\max}(\bSig_0) \leq c_0^{-1}$ for some constant $c_0>0$, then it holds that
\begin{align*}
&\|\wt{\bB} - \bI_{k}\|_2 \leq c_1\|\bOmg - \bOmg_0\|_2 = O(a_{n}),
\end{align*}
where $\wt{\bB}$ is given in (\ref{eq: proxy-knockoff}) and $c_1>0$ is some  uniform constant independent of set $\mathcal{S}$.
\end{lemma}

\noindent \textit{Proof}. We use $C$ to denote some generic positive constant whose value may change from line to line.  First note that
\begin{align}\label{eq: e025}
(\bB^{\sbOmg}_\mathcal{S})^T\bB^{\sbOmg}_\mathcal{S} - \bB_{0,\mathcal{S}}^T\bB_{0,\mathcal{S}} = - \big(\diag(\bs)(\bOmg-\bOmg_0)\diag(\bs)\big)_{\mathcal{S},\mathcal{S}}.
\end{align}
Further, since $\Sig_0 - 2^{-1}\diag(\bs) $ is positive definite it follows that $\|\bs\|_\infty \leq 2\Lambda_{\max}(\Sig_0) \leq 2c_0^{-1}$. Thus it holds that
\begin{align*}
\|(\bB^{\sbOmg}_\mathcal{S})^T\bB^{\sbOmg}_\mathcal{S} - \bB_{0,\mathcal{S}}^T\bB_{0,\mathcal{S}}\|_2 \leq C\|(\bOmg-\bOmg_0)_{\mathcal{S},\mathcal{S}}\|_2 \leq C \|\bOmg-\bOmg_0\|_2=O(a_{n}).
\end{align*}
For $n$ large enough, by the triangle inequality we have
\begin{align*}
 \Lambda_{\min}((\bB^{\sbOmg}_\mathcal{S})^T\bB^{\sbOmg}_\mathcal{S} ) & \geq \Lambda_{\min}(\bB_{0,\mathcal{S}}^T\bB_{0,\mathcal{S}} ) + \Lambda_{\min}((\bB^{\sbOmg}_\mathcal{S})^T\bB^{\sbOmg}_\mathcal{S}-\bB_{0,\mathcal{S}}^T\bB_{0,\mathcal{S}} ) \\
& \geq \Lambda_{\min}(\bB_{0}^T\bB_{0} ) - O(a_{n}) \\
& = \Lambda_{\min}\Big(2 \diag(\bs) - \diag(\bs)\bOmg_0\diag(\bs) \Big) - O(a_{n})\\
&\geq c_0/2.
\end{align*}
In addition, $ \Lambda_{\min}((\bB_{0,\mathcal{S}})^T\bB_{0,\mathcal{S}}) \geq \Lambda_{\min}((\bB_{0})^T\bB_{0}) = \Lambda_{\min}\Big(2 \diag(\bs) - \diag(\bs)\bOmg_0\diag(\bs) \Big) \geq c_0/2$. The above two inequalities together with Lemma 2.2 in \cite{schmitt1992perturbation} entail that
\begin{align}
\nonumber\left\|\big((\bB^{\sbOmg}_\mathcal{S})^T\bB^{\sbOmg}_\mathcal{S} \big)^{1/2} - \big((\bB_{0,\mathcal{S}})^T\bB_{0,\mathcal{S}}\big)^{1/2} \right\|_2 & \leq (\sqrt{c_0/2} + \sqrt{c_0/2})^{-1}\|(\bB^{\sbOmg}_\mathcal{S})^T\bB^{\sbOmg}_\mathcal{S} - (\bB_{0,\mathcal{S}})^T\bB_{0,\mathcal{S}} \|_2 \\
& \leq C \|\bOmg - \bOmg_0\|_2= O(a_{n}),
\end{align}
where the last step is because of \eqref{eq: e025}. Thus it follows that
\begin{align}\label{eq: e054}
\nonumber \|\wt{\bB} - \bI_k\|_2 &\leq  \left\|\big((\bB^{\sbOmg}_\mathcal{S})^T\bB^{\sbOmg}_\mathcal{S} \big)^{1/2} - \big((\bB_{0,\mathcal{S}})^T\bB_{0,\mathcal{S}}\big)^{1/2} \right\|_2 \|\big((\bB_{0,\mathcal{S}})^T\bB_{0,\mathcal{S}}\big)^{-1/2}\|_2 \\
& \leq C \|\bOmg - \bOmg_0\|_2 \Lambda_{\min}^{-1/2}( (\bB_{0})^T\bB_{0}  ) \leq C\|\bOmg - \bOmg_0\|_2 ,
\end{align}
where the last step comes from assumption $\Lambda_{\min}(\bB_{0}^T\bB_0)= \Lambda_{\min}(2\diag(\bs) - \diag(\bs)\bOmg_0\diag(\bs)) \\ \geq c_0$.  This concludes the proof of Lemma \ref{lemma: 1}.

\subsection{Proof of Theorem \ref{thm1}} \label{appA.4}

We now proceed to prove Theorem \ref{thm1} with the aid of Lemma \ref{lemma: 1} in Section \ref{appA.3}. We use the same notation  as in the proof of Proposition \ref{prop1} and use $C>0$ to denote a generic constant whose value may change from line to line.

We start with proving \eqref{eq: e042}. By Condition \ref{cond0}, we have
\begin{align}\label{eq: e041}
|{\rm FDR}(\bH^{\sbOmg}, \mathcal{S}) - {\rm FDR}(\bH^{\sbOmg_0}, \mathcal{S})| \leq L\|\bH^{\sbOmg} - \bH^{\sbOmg_0}\|_F,
\end{align}
where the constant $L$ is uniform over all $\|\bOmg - \bOmg_{0}\|\leq C_2a_n$ and $|\mathcal{S}|\leq K_n$.  Denote by $k = |\mathcal{S}|$. By the definition of $\bH^{\sbOmg}$, it holds that
\begin{align*}
\bH^{\sbOmg} - \bH^{\sbOmg_0} = \bH^{\sbOmg} - \bI_{p+k} = \left(\begin{array}{cc}
\bzero & \bC_{\mathcal S}^{\sbOmg}- \bC_{0,\mathcal{S}}\wt{\bB} \\
\bzero & \wt{\bB} -\bI_k
\end{array}\right).
\end{align*}
By the definition and matrix norm inequality, we deduce
\begin{align*}
\|\bH^{\sbOmg} - \bH^{\sbOmg_0} \|_F  & = \|\bC_{\mathcal S}^{\sbOmg}- \bC_{0,\mathcal{S}}\wt{\bB} \|_F + \|\wt{\bB} - \bI_k\|_F\\
&\leq  \sqrt{k}\|\bC_{\mathcal S}^{\sbOmg}- \bC_{0,\mathcal{S}}\wt{\bB} \|_2 + \sqrt{k} \|\wt{\bB} - \bI_k\|_2 \\
&\leq \sqrt{K_n}\Big(\|\bC_{\mathcal S}^{\sbOmg}- \bC_{0,\mathcal{S}}\|_2 + \|\bC_{0,\mathcal{S}}(\wt{\bB} -\bI_k) \|_2 + \|\wt{\bB} - \bI_k\|_2\Big)\\
& \leq \sqrt{K_n}\Big( \big(1+\|\bC_{0,\mathcal S}\|_2\big)\|\wt{\bB} - \bI_k\|_2 + \|\bC^{\sbOmg}- \bC_{0}\|_2\Big). 
\end{align*}
Since $\bSig_0 - 2^{-1}\diag(\bs)$ is positive definite, it follows that $s_j \leq 2\Lambda_{\max}(\Sig_0) \leq 2/\Lambda_{\min}(\bOmg_0)\leq C$. Thus  $\|\bC_{0,\mathcal S}\|_2 \leq \|\bC_{0}\|_2 = \|\bI - \bOmg_{0}\diag(\bs)\|_2 \leq 1 + \| \bOmg_{0}\|_2\|\diag(\bs)\|_2 \leq C$. This along with $\|\bC^{\sbOmg}- \bC_{0}\|_2 = \|(\bOmg-\bOmg_{0})\diag(\bs)\|_2\leq Ca_n$ and
Lemma \ref{lemma: 1} entails that $\|\bH^{\sbOmg} - \bH^{\sbOmg_0} \|_F$ can be further bounded as
\[
\|\bH^{\sbOmg} - \bH^{\sbOmg_0} \|_F \leq C\sqrt{K_n}a_n.
\]
Combining the above result with \eqref{eq: e041} leads to
\begin{align}
\sup_{\{|\mathcal{S}|\leq K_n, \|\sbOmg - \sbOmg_0\|\leq Ca_n\}}|{\rm FDR}(\bH^{\sbOmg}, \mathcal{S}) - {\rm FDR}(\bH^{\sbOmg_0}, \mathcal{S})| \leq O(\sqrt{K_n}a_n),
\end{align}
which completes the proof of \eqref{eq: e042}.

We next establish the FDR control for RANK. By Condition \ref{cond5}, the event $\mathcal{E}_0 = \{\|\wh{\bOmg} - \bOmg_{0}\|_2\leq C_2a_n\}$ occurs with probability at least $1-O(p^{-c_1})$. Since $\wh{\bOmg}$ and $\wt{\mathcal{S}}$ are estimated from independent subsample $(\bX^{(1)}, \by^{(1)})$, it follows from \eqref{eq: e042} that
\begin{align}
\nonumber& \left|{\mb E}\big[{\rm FDP}_n(\bX_{\wt{\mathcal{S}}}^{(2)}, \hbX_{\wt{\mathcal{S}}})\big| \mathcal{E}_0\big] - {\mb E}\big[{\rm FDP}_n(\bX_{\wt{\mathcal{S}}}^{(2)}, \tbX_{0,\wt{\mathcal{S}}})\big| \mathcal{E}_0\big]\right| \\
\nonumber \leq &\sup_{|\mathcal{S}|\leq K_n,\, \|\sbOmg - \sbOmg_0\|\leq C_2a_n}\left|{\mb E}\big[{\rm FDP}_n(\bX_{\mathcal{S}}^{(2)}, \tbX_{\mathcal{S}}^{\sbOmg})\big| \mathcal{E}_0\big] - {\mb E}\big[{\rm FDP}_n(\bX_{\mathcal{S}}^{(2)}, \tbX_{0,\mathcal{S}})\big| \mathcal{E}_0\big]\right| \\
\nonumber= &\sup_{|\mathcal{S}|\leq K_n,\, \|\sbOmg - \sbOmg_0\|\leq C_2a_n}\left|{\mb E}\big[{\rm FDP}_n(\bX_{\mathcal{S}}^{(2)}, \tbX_{\mathcal{S}}^{\sbOmg})\big] - {\mb E}\big[{\rm FDP}_n(\bX_{\mathcal{S}}^{(2)}, \tbX_{0,\mathcal{S}})\big]\right|\\
\nonumber= &\sup_{|\mathcal{S}|\leq K_n,\, \|\sbOmg - \sbOmg_0\|\leq C_2a_n}\left|{\rm FDR}(\bH^{\sbOmg}, \mathcal{S}) - {\rm FDR}(\bH^{\sbOmg_0}, \mathcal{S})\right|\\
 \leq & O(\sqrt{K_n}a_n).\label{eq: e052}
\end{align}
Now note that by the property of conditional expectation, we have
\begin{align*}
 &{\rm FDR}_n(\wh{\bOmg}, \wt{\mathcal{S}}) - {\rm FDR}_n(\bOmg_0, \wt{\mathcal{S}}) \\
 & =
\left( {\mb E}\big[{\rm FDP}_n(\bX_{\wt{\mathcal{S}}}^{(2)}, \hbX_{\wt{\mathcal{S}}})\big| \mathcal{E}_0\big] - {\mb E}\big[{\rm FDP}_n(\bX_{\wt{\mathcal{S}}}^{(2)}, \tbX_{0,\wt{\mathcal{S}}})\big| \mathcal{E}_0\big]\right) {\mb P}(\mathcal{E}_0) \\
& \quad +  \left({\mb E}\big[{\rm FDP}_n(\bX_{\wt{\mathcal{S}}}^{(2)}, \hbX_{\wt{\mathcal{S}}})\big| \mathcal{E}_0^c\big] - {\mb E}\big[{\rm FDP}_n(\bX_{\wt{\mathcal{S}}}^{(2)}, \tbX_{0,\wt{\mathcal{S}}})\big| \mathcal{E}_0^c\big] \right){\mb P}(\mathcal{E}_0^c) \\
& \equiv I_1 + I_2.
\end{align*}
Let us first consider term $I_1$. By \eqref{eq: e052}, it holds that
\begin{align*}
|I_1| \leq \left| {\mb E}\big[{\rm FDP}_n(\bX_{\wt{\mathcal{S}}}^{(2)}, \hbX_{\wt{\mathcal{S}}})\big| \mathcal{E}_0\big] - {\mb E}\big[{\rm FDP}_n(\bX_{\wt{\mathcal{S}}}^{(2)}, \tbX_{0,\wt{\mathcal{S}}})\big| \mathcal{E}_0\big]\right| \leq O(\sqrt{K_n}a_n).
\end{align*}
We next consider term $I_2$. Since FDP is always bounded between 0 and 1, we have
\begin{align*}
|I_2| \leq 2{\mb P}(\mathcal{E}_0^c) \leq O(p^{-c_1}).
\end{align*}
Combining the above two results yields
\begin{align*}
 \left|{\rm FDR}_n(\wh{\bOmg}, \wt{\mathcal{S}}) - {\rm FDR}_n(\bOmg_0, \wt{\mathcal{S}}) \right|  \leq O(\sqrt{K_n}a_n) + O(p^{-c_1}).
\end{align*}
This together with the result of ${\rm FDR}_n(\bOmg_0, \wt{\mathcal{S}}) \leq q$ mentioned in the proof of Proposition \ref{prop1} in Section \ref{appA.2} completes 
the proof of Theorem \ref{thm1}.

\subsection{Proof of Theorem \ref{theo2}} \label{appA.5}

In this proof, we will drop the superscripts in $(\bX^{(2)}, \by^{(2)})$ whenever there is no confusion.
By the definition of power, for any given precision matrix $\bOmg$ and reduced model $\mathcal{S}$ the power can be written as
\begin{align*}
\text{Power}(\bOmg, \mathcal{S}) = {\mb E}[f(\bX_{\mathcal{S}}, \tbX^{\sbOmg}_{\mathcal{S}}, \by)],
\end{align*}
where $f$ is some function describing how the empirical power depends on the data. Note that $f(\bX_{\mathcal{S}}, \tbX^{\sbOmg}_{\mathcal{S}}, \by)$ is a stochastic process indexed by $\bOmg$, and we care about the mean of this process. Our main idea is to construct another stochastic process indexed by $\bOmg$ which has the same mean but possibly different
distribution. Then by studying the mean of this new stochastic process, we can prove the desired result.

We next provide more technical details of the proof. The proxy process is defined as
\begin{align}\label{eq: proxy-knockoff-power}
& \breve{\bX}^{\sbOmg}_{\mathcal{S}}  = \bX\bC_{\mathcal{S}}^{\sbOmg} + \bZ\bB_{0,\mathcal{S}}(\bB_{0,\mathcal{S}}^T\bB_{0,\mathcal{S}})^{-1/2}\Big(\big(\bB_{\mathcal{S}}^{\sbOmg}\big)^T\bB_{\mathcal{S}}\Big)^{1/2},
\end{align}
where $\bC^{\sbOmg}_{\mathcal S}$ is the submatrix of $\bC^{\sbOmg}=\bI_p- \bOmg\diag\{\bs\}$, $\bB^{\sbOmg}_{\mathcal{S}}$ is the submatrix of $\bB^{\sbOmg} = \big(\diag(\bs) - \diag(\bs)\bOmg\diag(\bs)\big)^{1/2}$, and $\bB_0 = \bB^{\sbOmg_0}$. It is easy to see that $\breve{\bX}^{\sbOmg}_{\mathcal{S}}$ and $
\widetilde{\bX}_{\mathcal{S}}^{\sbOmg}$ defined using \eqref{eq: xtilde-distr} have the same distribution. Since $\bZ$ is independent of $(\bX, \by)$, we can further conclude that $(\bX_{\mathcal{S}}, \tbX_{\mathcal{S}}^{\sbOmg}, \by)$ and $(\bX_{\mathcal{S}}, \breve{\bX}_{\mathcal{S}}^{\sbOmg}, \by)$ have the same joint distribution for each given $\bOmg$ and $\mathcal S$. Thus the power function can be further written as
\begin{align*}
\text{Power}(\bOmg, \mathcal{S}) = {\mb E}[f(\bX_{\mathcal{S}}, \tbX^{\sbOmg}_{\mathcal{S}}, \by)] = {\mb E}[f(\bX_{\mathcal{S}}, \breve\bX^{\sbOmg}_{\mathcal{S}}, \by)].
\end{align*}
Therefore, we only need to study the power of the knockoffs procedure based on the pseudo data $(\bX_{\mathcal{S}}, \breve{\bX}_{\mathcal{S}}^{\sbOmg}, \by)$.

To simplify the technical presentation, we will slightly abuse the notation and still use $\widehat{\bbeta}=  \widehat{\bbeta}(\lambda)  = \widehat{\bbeta}(\lambda; \bOmg, \mathcal{S})$ to represent the Lasso solution based on pseudo data $(\bX_{\mathcal{S}}, \breve{\bX}_{\mathcal{S}}^{\sbOmg}, \by)$. We will use $c$ and $C$ to denote some generic positive constants whose values may change from line to line.
Define
\begin{eqnarray} \label{eqn:def-Gtilde}
	\widetilde{\bG}=\dfrac{1}{n}\wt{\bX}_{\rm KO}^{T}\wt{\bX}_{\rm KO}\in \mb{R}^{(2p)\times (2p)} 
	\quad \text{ and }  \quad  \widetilde{\brho}=\dfrac{1}{n}\wt{\bX}_{\rm KO}^{T}\by \in \mb{R}^{2p}
\end{eqnarray}
with $\wt{\bX}_{\rm KO} = [\bX, \breve{\bX}^{\sbOmg}] \in \mathbb{R}^{n\times(2p)}$ the augmented design matrix. For any given set $\mathcal{S}\subset \{1,\cdots, p\}$ with $k = |\mathcal{S}|$, $(2p)\times (2p)$ matrix $\bA$, and $(2p)$-vector $\ba$, we will abuse the notation and denote by $\bA_{\mathcal{S},\mathcal{S}}\in \mb{R}^{(2k) \times (2k)}$ the principal submatrix formed by columns and rows in set $\{j: j\in \mathcal{S} \text{ or } j-p \in \mathcal{S}\}$ and $\ba_{\mathcal{S}} \in \mathbb{R}^{2k}$ the subvector formed by components in set $\{j: j\in \mathcal{S} \text{ or } j-p \in \mathcal{S}\}$. For any $p\times p$ matrix $\bB$ (or $p$-vector $\bb$), we define $\bB_{\mathcal S}$ (or $\bb_{\mathcal{S}}$) in the same way meaning that columns (or components) in set $\mathcal{S}$ will be taken to form the submatrix (or subvector).

With the above notation, note that the Lasso solution $\hbbeta = (\hbeta_1,\cdots, \hbeta_{2p})^T = \hbbeta(\lambda; \bOmg, \mathcal S)$ restricted to variables in $\mathcal{S}$ can be  obtained by setting $\widehat{\beta}_j = 0$ for $j \in  \{1\leq j \leq 2p: j\not\in \mathcal{S} \text{ and } j-p\not\in \mathcal{S}\} $ and minimizing the following objective function
\begin{eqnarray}\label{lasso-obj-ko}
\widehat{\bbeta}_{\mathcal{S}}={\arg\min}_{\bb\in \mathbb{R}^{2k}}\left\{\dfrac{1}{2}\bb^T\widetilde{\bG}_{\mathcal{S},\mathcal{S}}\bb-\widetilde{\brho}_{\mathcal{S}}^T\bb+\lambda\|\bb\|_1\right\}.
\end{eqnarray}
By Proposition \ref{pro1} in Section \ref{appA.6}, with probability at least $1-O(p^{-c})$ it holds that
\begin{align}
& \sup_{\sbOmg\in \calA,\, |\mathcal S|\leq K_n} \|\widehat{\bbeta}(\lambda; \bOmg, \mathcal S) - \bbeta_{\mb T} \|_2= O(\sqrt{s}\lambda), \label{eq: e056}\\
& \sup_{\sbOmg\in \calA,\, |\mathcal S|\leq K_n} \|\widehat{\bbeta}(\lambda; \bOmg, \mathcal S)  - \bbeta_{\mb T} \|_1= O(s\lambda)\label{eq: e057},
\end{align}
where $\lambda = C\sqrt{(\log p)/n} $ with $C > 0$ some constant.

Denote by $W_j^{\sbOmg, \mathcal S}$ the LCD based on the above $\widehat{\bbeta}(\lambda; \bOmg, \mathcal S)$.
Recall that by assumption, there are no ties in the magnitude of nonzero $W_j^{\sbOmg, \mathcal S}$'s and no ties in the nonzero components of the Lasso solution with asymptotic probability one. Let $|W_{(1)}^{\sbOmg, \mathcal S}| \geq \cdots \geq |W_{(p)}^{\sbOmg,\mathcal S}|$ be the ordered knockoff statistics according to magnitude. Denote by $j^*$ the index such that $|W_{(j^*)}^{\sbOmg, \mathcal S}| = T$. Then by the definition of $T$, it holds that $-T <W_{(j^*+1)}^{\sbOmg, \mathcal S} \leq 0$. We next analyze  the two cases of $W_{(j^*+1)}^{\sbOmg, \mathcal S} = 0$ and $-T<W_{(j^*+1)}^{\sbOmg, \mathcal S} < 0$ separately.

\textbf{Case 1.} For the case of $W_{(j^*+1)}^{\sbOmg, \mathcal S} = 0$, we have $ W_{(j^*+1)}^{\sbOmg, \mathcal S}= \cdots = W_{(p)}^{\sbOmg, \mathcal S}=0$.   It follows from the definition of $T$ that $\wh{\mathcal S}^{\sbOmg} = \supp(\bW^{\sbOmg, \mathcal S}) $ and thus
\begin{align}\label{eq: e055}
\wh{\mathcal S}^{\sbOmg} \supset \{1,\cdots, p\}\setminus\mathcal{S}_1^{\sbOmg},
\end{align}
where $\mathcal{S}_1^{\sbOmg} =\{1\leq j\leq p: \hbeta_j(\lambda;\bOmg, \mathcal S) = 0\}$.
Meanwhile, note that in view of \eqref{eq: e057} we have with probability at least $1-O(p^{-c})$,
\begin{align*}
O(s\lambda) &\geq \sup_{\sbOmg\in \calA,\, |\mathcal S|\leq K_n}\|\hbbeta(\lambda; \bOmg, \mathcal S) - \bbeta_0\|_1 \geq \sup_{\sbOmg\in \calA,\, |\mathcal S|\leq K_n} \sum_{j\in \mathcal{S}_1^{\sbOmg}\cap \mathcal{S}_0}|\hbeta_{j}(\lambda; \bOmg, \mathcal S) - \beta_{0,j}|\\
& =  \sum_{j\in \mathcal{S}_1^{\sbOmg}\cap \mathcal{S}_0}|\beta_{0,j}| \geq |\mathcal{S}_1^{\sbOmg}\cap \mathcal{S}_0| \min_{j\in \mathcal{S}_0}|\beta_{0,j}|.
\end{align*}
By Condition \ref{cond3} and $\lambda = O(\sqrt{(\log p)/n})$, we can further derive from the above inequality that
\[
|\mathcal{S}_1^{\sbOmg}\cap \mathcal{S}_0| = o(s),
\]
which together with $|\mathcal{S}_0| = s$ entails that
\[
|\big(\{1,\cdots,p\}\setminus\mathcal{S}_1^{\sbOmg}\big)\cap \mathcal{S}_0| \geq [1-o(1)]s.
\]
Combining this with \eqref{eq: e055} leads to
\[
\Big|\wh{\mathcal S}^{\sbOmg} \cap \mathcal{S}_0\Big| \geq \Big |\big(\{1,\cdots,p\}\setminus\mathcal{S}_1^{\sbOmg}\big)\cap \mathcal{S}_0\Big| \geq [1-o(1)]s.
\]
Thus with asymptotic probability one, it holds uniformly over all $\bOmg \in \calA$ and $|\mathcal S| \leq K_n$ that
\[
\dfrac{\left|\wh{\mathcal S}^{\sbOmg} \cap \mathcal{S}_0\right| }{s}
 \geq 1-o(1).
\]

\textbf{Case 2.} We consider the case of $-T<W_{(j^*+1)}^{\sbOmg, \mathcal S} < 0$. By the definitions of $T$ and $j^*$, we have
\begin{align}\label{eq: e051}
\frac{|\{j: W_j^{\sbOmg, \mathcal S} \leq -T\}| + 2}{|\{j: W_j^{\sbOmg, \mathcal S} \geq T\}|} > q
\end{align}
since otherwise we would reduce $T$ to $|W_{(j^*+1)}^{\sbOmg, \mathcal S}|$ to get the new smaller threshold with the criterion still satisfied.
We next bound $T$ using the results in Lemma \ref{lem: cond-Shat} in Section \ref{appB.4} of Supplementary Material. Observe that \eqref{eq: e051} and Lemma \ref{lem: cond-Shat} lead to
$ |\{j : W_j^{\sbOmg,  \mathcal S}\leqslant - T\} | >
q |\{j : W_j^{\sbOmg,  \mathcal S}\geqslant T\} | -2 \geq
q c s -2$
with asymptotic probability one. Moreover, when $W_{j}^{\sbOmg, \mathcal S} \leq - T$ we have $|\widehat{\beta}_{j}(\lambda;\bOmg, \mathcal S)| - |\widehat{\beta}_{j + p}(\lambda; \bOmg, \mathcal S)| \leq - T $ and thus  $|\widehat{\beta}_{j+p}(\lambda; \bOmg, \mathcal S)| \geq T $.
Using \eqref{eq: e057}, we obtain
\begin{align*}
O(s\lambda)  =
\|\widehat{\bbeta}(\lambda; \bOmg, \mathcal S)-\bbeta_{\mathbb{T}}\|_1
\geq \sum _{j :\, W_j^{\sbOmg, \mathcal S}\leqslant - T} | \widehat{\beta}_{j+p}(\lambda; \bOmg, \mathcal S) |
\geq T |\{j : W_j^{\sbOmg, \mathcal S}\leqslant - T\}| .
\end{align*}
Combining these results yields
$O(s\lambda) \geq T (q c s-1)$ and thus it holds that
\begin{equation}
T \leq O(\lambda).
\end{equation}

We now proceed to prove the theorem by showing that Type II error is small.
By Condition \ref{cond3} of $\sqrt{n/(\log p)}\min\limits_{j\in \mathcal{S}_0}|\beta_{0,j}|\rightarrow \infty$ and assumption $\lambda= C\sqrt{(\log p)/n}$,  there exists some $\kappa_n \rightarrow \infty$ such that $\min\limits_{j\in \mathcal{S}_0}|\beta_{0,j}|\geq \kappa_n\lambda$ as $n\rightarrow \infty$.
In light of \eqref{eq: e057}, we derive
\begin{align*}
O(s\lambda)
& = \|\widehat{\bbeta}(\lambda; \bOmg, \mathcal S)-\bbeta_{\mathbb{T}}\|_1
= \sum_{j=1}^{p} (| \widehat{\beta}_{j}(\lambda; \bOmg, \mathcal S) - \beta_{j} |+ | \widehat{\beta}_{j+p}(\lambda; \bOmg, \mathcal S) | )\\
& \geq \sum_{j \in \mathcal{S}_0 \cap (\widehat{\mathcal S}^{\sbOmg})^{c} } (| \widehat{\beta}_{j}(\lambda; \bOmg, \mathcal S) - \beta_{j} | + | \widehat{\beta}_{j+p}(\lambda; \bOmg, \mathcal S) | )\\
& \geq \sum_{j \in \mathcal{S}_0 \cap (\widehat{\mathcal S}^{\sbOmg})^{c} } (| \widehat{\beta}_{j}(\lambda; \bOmg, \mathcal S) - \beta_{j} | + | \widehat{\beta}_{j}(\lambda; \bOmg, \mathcal S) |  - T)
\end{align*}
since $| \widehat{\beta}_{j+p}(\lambda; \bOmg, \mathcal S)| \geq | \widehat{\beta}_{j}(\lambda; \bOmg, \mathcal S) | - T $ when $j \in (\widehat{\mathcal S}^{\sbOmg}) ^{c}$.
Using the triangle inequality and since $ |\beta_{j}| \geq \lambda\kappa_{n} $ when $j \in \mathcal{S}_0 $, we can conclude that
$$ O(s\lambda)
\geq \sum_{j \in \mathcal{S}_0 \cap (\widehat{\mathcal S}^{\sbOmg})^{c} } (| \beta_{j} | - T)
\geq (\lambda\kappa_{n} - T) \left|\{j\in \mathcal{S}_0\cap (\widehat{\mathcal S}^{\sbOmg})^{c}\}\right|. $$
Then it follows that
\begin{align*}
 \dfrac{ |\{j\in \mathcal{S}_0\cap (\widehat{\mathcal S}^{\sbOmg})^c\}|}{s}
& = 1 - \dfrac{ |\{j\in \mathcal{S}_0\cap (\widehat{\mathcal S}^{\sbOmg})^{c}\}| }{s}
\geq 1 - O \left( \dfrac{\lambda}{\lambda\kappa_{n} - T} \right) \\
& = 1 - O(\kappa_{n} ^{-1}) = 1-o(1)
\end{align*}
uniformly over all $\bOmg \in \calA$ and $|\mathcal S| \leq K_n$.

Combining the above two scenarios, we have shown that with asymptotic probability one, uniformly over all $\bOmg \in \calA$ and $|\mathcal S| \leq K_n$ it holds that
\begin{eqnarray}\label{eq: e058}
 \dfrac{ |\{j\in \mathcal{S}_0\cap (\widehat{\mathcal S}^{\sbOmg})^c\}|}{s}
\geq 1-o(1).
\end{eqnarray}
This along with the assumption $\mathbb{P}\{\wh{\bOmg} \in \calA\} =  1-O(p^{-c_2})$ in Condition \ref{cond5} gives
\begin{align*}
{\rm Power}(\wh{\bOmg}, \wt{\mathcal S}) &= \mb{E}\Big[\dfrac{ |\{j\in \mathcal{S}_0\cap (\widehat{\mathcal S}^{\wh{\sbOmg}})^c\}|}{s}\Big] \\
& \geq \mb{E}\Big[\dfrac{ |\{j\in \mathcal{S}_0\cap (\widehat{\mathcal S}^{\wh{\sbOmg}})^c\}|}{s}\Big| \wh{\bOmg}\in \calA\Big] \mathbb{P}\{\wh{\bOmg}\in \calA\} \\
& \geq [1-o(1)][1-O(p^{-c_2})] = 1 -o(1),
\end{align*}
where the second to the last step is due to \eqref{eq: e058}. This concludes the proof of Theorem \ref{theo2}.

\subsection{Proposition \ref{pro1} and its proof} \label{appA.6}

\begin{proposition}\label{pro1}
Assume that Conditions \ref{cond4} and \ref{cond2} hold, the smallest  eigenvalue of $2\diag(\bs) - \diag(\bs)\bOmg_0\diag(\bs)$ is positive and bounded away from $0$, and $\lambda =  C_3\{(\log p)/n\}^{1/2}$ with $C_3>0$ some constant. Let $\bbeta_{\mathbb{T}}=(\bbeta_0^T, 0,\cdots,0)^T \in \mathbb{R}^{2p}$ be the expanded vector of true regression coefficient vector. If $[(L_p + L_p^{\prime})^{1/2}+K_n^{1/2}] a_n = o(1)$ and $s \{a_n + L_{p}' [(\log p)/n]^{1/2} + [K_nL_p'(\log p)/n]^{1/2} \} = o(1)$, then with probability at least $1-O(p^{-c_3})$,
\[
		\sup_{\sbOmg \in \calA,\, |\mathcal{S}|\leq K_n}\|\widehat{\bbeta}(
		\lambda; \bOmg, \mathcal{S})-\bbeta_{\mathbb{T}}\|_1=O(s\lambda) \text{ and }
		\sup_{\sbOmg \in \calA,\, |\mathcal{S}|\leq K_n}\|\widehat{\bbeta}(\lambda; \bOmg, \mathcal{S})-\bbeta_{\mathbb{T}}\|_2=O(s^{1/2}\lambda),
\]
where $\widehat{\bbeta}(\lambda; \bOmg, \mathcal{S})$ is defined in the proof of Theorem \ref{theo2} in Section \ref{appA.5} and $c_3 > 0$ is some constant.
\end{proposition}

\noindent \textit{Proof}. We adopt the same notation as used in the proof of Theorem \ref{theo2} in Section \ref{appA.5}. Let us introduce some key events which will be used in the technical analysis. Define
		\begin{eqnarray}\label{def: event3}
		\mathcal{E}_3&=&\Big\{\sup_{\|\sbOmg - \sbOmg_0\|\leq C_2a_n, \, |\mathcal{S}|\leq K_n}\|\widetilde{\brho}_{\mathcal{S}}-\widetilde{\bG}_{\mathcal{S},\mathcal{S}}\bbeta_{\mathbb{T}, \mathcal{S}}\|_{\infty}\le \lambda_0\Big\},\\
		\mathcal{E}_4&=&\Big\{\sup_{\|\sbOmg - \sbOmg_0\|_2\leq C_2a_n, \, |\mathcal{S}|\leq K_n}\|\wt{\bG}_{\mathcal{S},\mathcal{S}} - \bG_{\mathcal{S},\mathcal{S}}\|_{\max}\le C_5a_{2n}\Big\},
		\end{eqnarray}
		where  $\lambda_0 = C_4\sqrt{(\log p)/n}$ and $a_{2n} = a_n +( L_{p}' + K_n) \sqrt{(\log p)/n} $ with $C_4, C_5 > 0$ some constants.
		Then by Lemmas \ref{lem: beta-est} and \ref{bound-lem} in Sections \ref{appB.2} and \ref{appB.5} of Supplementary Material,
			\begin{align}\label{eq: event3-prob}
			{\mb P}( \mathcal{E}_3 \cap \mathcal{E}_4) = 1 - O(p^{-c_0})
			\end{align}
for some constant $c_0>0$. Hereafter we will condition on the event $ \mathcal{E}_3 \cap \mathcal{E}_4$.
		
		Since $\wh{\bbeta}_{\mathcal{S}}$ is the minimizer of the objective function  in (\ref{lasso-obj-ko}), we have
		\begin{eqnarray*}
			\dfrac{1}{2}\widehat{\bbeta}_{\mathcal{S}}^T\widetilde{\bG}_{\mathcal{S}, \mathcal{S}}\widehat{\bbeta}_{\mathcal{S}}-\widetilde{\brho}_{\mathcal{S}}^T\widehat{\bbeta}_{\mathcal{S}}+\lambda\|\widehat{\bbeta}_\mathcal{S}\|_1\le \dfrac{1}{2}\bbeta_{\mathbb{T},\mathcal{S}}^{T}\widetilde{\bG}_{\mathcal{S}, \mathcal{S}}\bbeta_{\mathbb{T},\mathcal{S}}-\widetilde{\brho}_{\mathcal S}^T\bbeta_{\mathbb{T},\mathcal{S}}+\lambda\|\bbeta_{\mathbb{T},\mathcal{S}}\|_1.
		\end{eqnarray*}
Some routine calculations lead to
		\begin{eqnarray}\nonumber
		&&\dfrac{1}{2}(\widehat{\bbeta}_\mathcal{S}-\bbeta_{\mathbb{T},\mathcal{S}})^T\widetilde{\bG}_{\mathcal{S},\mathcal{S}}(\widehat{\bbeta}_\mathcal{S}-\bbeta_{\mathbb{T},\mathcal{S}})+\lambda\|\widehat{\bbeta}_\mathcal{S}\|_1\\\nonumber
		&\le& -\bbeta_{\mathbb{T},\mathcal{S}}^{T}\widetilde{\bG}_{\mathcal{S},\mathcal{S}}\widehat{\bbeta}_{\mathcal{S}}+\bbeta_{\mathbb{T},\mathcal{S}}^{T}\widetilde{\bG}_{\mathcal{S},\mathcal{S}}\bbeta_{\mathbb{T},\mathcal{S}}+\widetilde{\brho}_{\mathcal{S}}^T(\widehat{\bbeta}_\mathcal{S}-\bbeta_{\mathbb{T},\mathcal{S}})+\lambda\|\bbeta_{\mathbb{T},\mathcal{S}}\|_1\\\nonumber
		&=&(\widetilde{\brho}_\mathcal{S}-\widetilde{\bG}_{\mathcal{S},\mathcal{S}}\bbeta_{\mathbb{T},\mathcal{S}})^T(\widehat{\bbeta}_\mathcal{S}-\bbeta_{\mathbb{T},\mathcal{S}})+\lambda\|\bbeta_{\mathbb{T},\mathcal{S}}\|_1\\\label{equ-S-1}
		&\le& \|\widehat{\bbeta}_{\mathcal{S}}-\bbeta_{\mathbb{T},\mathcal{S}}\|_1\|\widetilde{\brho}_{\mathcal{S}}-\widetilde{\bG}_{\mathcal{S},\mathcal{S}}\bbeta_{\mathbb{T},\mathcal{S}}\|_{\infty}+\lambda\|\bbeta_{\mathbb{T},\mathcal{S}}\|_1.
		\end{eqnarray}
		Let $\widehat{\bdelta}=\widehat{\bbeta}-\bbeta_{\mathbb{T}}$. Then we can simplify (\ref{equ-S-1}) as
		\begin{eqnarray}\label{equ-S-2}
		\dfrac{1}{2}\widehat{\bdelta}_\mathcal{S}^T\widetilde{\bG}_{\mathcal{S},\mathcal{S}}\widehat{\bdelta}_\mathcal{S}+\lambda\|\widehat{\bbeta}_{\mathcal{S}}\|_1&\le& \|\widehat{\bdelta}_\mathcal{S}\|_1\|\widetilde{\brho}_{\mathcal{S}}-\widetilde{\bG}_{\mathcal{S},\mathcal{S}}\bbeta_{\mathbb{T},\mathcal{S}}\|_{\infty}+\lambda\|\bbeta_{\mathbb{T},\mathcal{S}}\|_1 \nonumber\\
		& \le&\lambda_0 \|\widehat{\bdelta}_\mathcal{S}\|_1 + \lambda\|\bbeta_{\mathbb{T},\mathcal{S}}\|_1.
		\end{eqnarray}
Observe that $\|\widehat{\bbeta}_\mathcal{S}\|_1=\|\widehat{\bbeta}_{\mathcal{S}_0}\|_1+\|\widehat{\bbeta}_{\mathcal{S}\setminus\mathcal{S}_0}\|_1$ and $\|\bbeta_{\mathbb{T},\mathcal{S}}\|_1=\|\bbeta_{\mathbb{T}, \mathcal{S}_0}\|_1+\|\bbeta_{\mathbb{T},\mathcal{S}\setminus\mathcal{S}_0}\|_1=\|\bbeta_{\mathbb{T}, \mathcal{S}_0}\|_1$ with $\mathcal{S}_0$ the support of true regression coefficient vector. Then it follows from $\|\widehat{\bbeta}_{\mathcal{S}_0}-\bbeta_{\mathbb{T}, \mathcal{S}_0}\|_1\geq\|\bbeta_{\mathbb{T}, \mathcal{S}_0}\|_1-\|\widehat{\bbeta}_{\mathcal{S}_0}\|_1$ that
		\begin{eqnarray*}
			\dfrac{1}{2}\widehat{\bdelta}_\mathcal{S}^T\widetilde{\bG}_{\mathcal{S},\mathcal{S}}\widehat{\bdelta}_\mathcal{S}+\lambda\|\widehat{\bbeta}_{\mathcal{S}\setminus\mathcal{S}_0}\|_1 \le \lambda_{0}\|\widehat{\bdelta}_\mathcal{S}\|_1+\lambda\|\widehat{\bbeta}_{\mathcal{S}_0}-\bbeta_{\mathbb{T}, \mathcal{S}_0}\|_1.
		\end{eqnarray*}
	
Denote by  $\|\widehat{\bdelta}_{\mathcal{S}_0}\|_1=\|\widehat{\bbeta}_{\mathcal{S}_0}-\bbeta_{\mathbb{T}, \mathcal{S}_0}\|_1$ and $\|\widehat{\bdelta}_{\mathcal{S}\setminus\mathcal{S}_0}\|_1=\|\widehat{\bbeta}_{\mathcal{S}\setminus \mathcal{S}_0}-\bbeta_{\mathbb{T}, \mathcal{S}\setminus\mathcal{S}_0}\|_1=\|\widehat{\bbeta}_{\mathcal{S}\setminus\mathcal{S}_0}\|_1$. Then we can further deduce
		\begin{eqnarray*}
			\dfrac{1}{2}\widehat{\bdelta}_\mathcal{S}^T\widetilde{\bG}_{\mathcal{S},\mathcal{S}}\widehat{\bdelta}_\mathcal{S}+\lambda\|\widehat{\bdelta}_{\mathcal{S}\setminus\mathcal{S}_0}\|_1\le \lambda_{0}\|\widehat{\bdelta}_\mathcal{S}\|_1+\lambda\|\widehat{\bdelta}_{\mathcal{S}_0}\|_1=\lambda_{0}\|\widehat{\bdelta}_{\mathcal{S}_0}\|_1+\lambda_{0}\|\widehat{\bdelta}_{\mathcal{S}\setminus\mathcal{S}_0}\|_1+\lambda\|\widehat{\bdelta}_{\mathcal{S}_0}\|_1;
		\end{eqnarray*}
		that is,
		\begin{eqnarray}\label{hatG-ineq}
		\dfrac{1}{2}\widehat{\bdelta}_\mathcal{S}^T\widetilde{\bG}_{\mathcal{S},\mathcal{S}}\widehat{\bdelta}_\mathcal{S}+(\lambda-\lambda_{0})\|\widehat{\bdelta}_{\mathcal{S}\setminus\mathcal{S}_0}\|_1 \le (\lambda+\lambda_{0})\|\widehat{\bdelta}_{\mathcal{S}_0}\|_1.
		\end{eqnarray}
		When $\lambda\geq 2\lambda_{0}$, it holds that
		\begin{eqnarray}\label{hatG-ineq-1}
		\dfrac{1}{2}\widehat{\bdelta}_\mathcal{S}^T\widetilde{\bG}_{\mathcal{S},\mathcal{S}}\widehat{\bdelta}_{\mathcal{S}}+\dfrac{\lambda}{2}\|\widehat{\bdelta}_{\mathcal{S}\setminus\mathcal{S}_0}\|_1 \le \dfrac{3\lambda}{2}\|\widehat{\bdelta}_{\mathcal{S}_0}\|_1.
		\end{eqnarray}
		Since  $\widehat{\bdelta}_\mathcal{S}^T\widetilde{\bG}_{\mathcal{S},\mathcal{S}}\widehat{\bdelta}_{\mathcal{S}}\geq0$,  we obtain the basic inequality
		\begin{align}\label{eq: e047}
		\|\widehat{\bdelta}_{\mathcal{S}\setminus\mathcal{S}_0}\|_1 \le 3\|\widehat{\bdelta}_{\mathcal{S}_0}\|_1
		\end{align}
		on event $\mathcal{E}_3$.
		It follows from \eqref{hatG-ineq-1} that
		\begin{eqnarray}\label{hatG-tildeG}
		\widehat{\bdelta}_\mathcal{S}^T{\bG}_{\mathcal{S},\mathcal{S}}\widehat{\bdelta}_\mathcal{S}
		\le 3\lambda\|\widehat{\bdelta}_{\mathcal{S}_0}\|_1+\widehat{\bdelta}_\mathcal{S}^T({\bG}_{\mathcal{S},\mathcal{S}}-\widetilde{\bG}_{\mathcal{S},\mathcal{S}})\widehat{\bdelta}_{\mathcal{S}}
		\end{eqnarray}
		with
		\[
		\bG = \left(\begin{array}{cc}\bSig_0 & \bSig_0-\diag(\bs) \\ \bSig_0-\diag(\bs)  &  \bSig_0 \end{array}\right).
		\]
With some matrix calculations, we can show that
		\begin{align*}
		\Lambda_{\min}(\bG) \geq C\Lambda_{\min}\{2\diag(\bs) -\diag(\bs)\bOmg_{0}\diag(\bs)\} \geq Cc_0,
		\end{align*}
		where the last step is by assumption.
		Thus the left hand side of (\ref{hatG-tildeG}) can be bounded from below by $c_0C\|\widehat{\bdelta}_\mathcal{S}\|_2^2$.
		
It remains to bound the right hand side of (\ref{hatG-tildeG}). For the first term, it follows from the Cauchy--Schwarz inequality and $2ab\le a^2/4+4b^2$ that
		\begin{eqnarray}\nonumber
		3\lambda\|\widehat{\bdelta}_{\mathcal{S}_0}\|_1&\le& 3\lambda\sqrt{s}\|\widehat{\bdelta}_{\mathcal{S}_0}\|_2\le 3\lambda\sqrt{s}\|\widehat{\bdelta}_\mathcal{S}\|_2 \\ \label{equ-46} &\le & 3\lambda\sqrt{s}
		\sqrt{\dfrac{\widehat{\bdelta}_\mathcal{S}^T\bG_{\mathcal{S},\mathcal{S}}\widehat{\bdelta}_\mathcal{S}}{cC_0}}
		\le C\lambda\sqrt{s}\|\wh{\bdelta}_\mathcal{S}\|_2, \label{eq: e049}
		\end{eqnarray}
		where the last step is because $\|\bG_{\mathcal{S},\mathcal{S}}\|_2 \leq 2\|\bSig_0\|_2 + 2\|\bSig_0 - \diag(\bs)\|_2 \leq 4\|\bSig_0\|_2 + \|\bs\|_\infty \leq C$ uniformly over all $\mathcal{S}$.
		For the last term $\widehat{\bdelta}_\mathcal{S}^T({\bG}_{\mathcal{S},\mathcal{S}}-\widetilde{\bG}_{\mathcal{S},\mathcal{S}})\widehat{\bdelta}_\mathcal{S}$ on the right hand of (\ref{hatG-tildeG}), by conditioning on event $\mathcal{E}_4$ and using the Cauchy--Schwarz inequality, the triangle inequality, and the basic inequality \eqref{eq: e047} we can obtain
		\begin{eqnarray}\nonumber
		\Big|\widehat{\bdelta}_\mathcal{S}^T({\bG}_{\mathcal{S},\mathcal{S}}-\widetilde{\bG}_{\mathcal{S},\mathcal{S}})\widehat{\bdelta}_\mathcal{S}\Big|&\le& \|{\bG}_{\mathcal{S},\mathcal{S}}-\widetilde{\bG}_{\mathcal{S},\mathcal{S}}\|_{\max}\|\widehat{\bdelta}_\mathcal{S}\|_1^2\\\nonumber
		&\le&\|{\bG}_{\mathcal{S},\mathcal{S}}-\widetilde{\bG}_{\mathcal{S},\mathcal{S}}\|_{\max}(\|\widehat{\bdelta}_{\mathcal{S}_0}\|_1+\|\widehat{\bdelta}_{\mathcal{S}\setminus\mathcal{S}_0}\|_1)^2\\\nonumber
		&\le&16\|{\bG}_{\mathcal{S},\mathcal{S}}-\widetilde{\bG}_{\mathcal{S},\mathcal{S}}\|_{\max}\|\widehat{\bdelta}_{\mathcal{S}_0}\|_1^2\\\nonumber
		&\le&16s\|{\bG}_{\mathcal{S},\mathcal{S}}-\widetilde{\bG}_{\mathcal{S},\mathcal{S}}\|_{\max}\|\widehat{\bdelta}_{\mathcal{S}_0}\|_2^2\le 16s\|{\bG}_{\mathcal{S},\mathcal{S}}-\widetilde{\bG}_{\mathcal{S},\mathcal{S}}\|_{\max}\|\widehat{\bdelta}_{\mathcal{S}}\|_2^2\\\nonumber
		&=&Csa_{2n}\|\widehat{\bdelta}_{\mathcal{S}}\|_2^2.
		\end{eqnarray}
		
Combining the above results, we can reduce inequality (\ref{hatG-tildeG}) to
		\begin{align*}
		c_0C\|\wh{\bdelta}_\mathcal{S}\|_2^2 \leq C\lambda\sqrt{s}\|\wh{\bdelta}_\mathcal{S}\|_2 + Csa_{2n}\|\wh{\bdelta}_\mathcal{S}\|_2^2.
		\end{align*}
		Since $sa_{2n}\rightarrow 0$, it holds for $n$ large enough that
		\begin{eqnarray*}
			\|\wh{\bdelta}_\mathcal{S}\|_2 = \|\widehat{\bbeta}_\mathcal{S}-\bbeta_{\mathbb{T},\mathcal{S}}\|_2=O(\sqrt{s}\lambda).
		\end{eqnarray*}
Further, by \eqref{eq: e049} we have
		\begin{align*}
		\|\widehat{\bbeta}_\mathcal{S}-\bbeta_{\mathbb{T},\mathcal{S}}\|_1=O(s\lambda).
		\end{align*}
Note that by definition, $\widehat{\bbeta}_{\mathcal{S}^c} = \bzero$ and $\bbeta_{\mathbb{T},\mathcal{S}^c}= \bzero$. Therefore, summarizing the above results completes the proof of Proposition \ref{pro1}.

\bibliographystyle{chicago}
\bibliography{references}


\newpage
\appendix
\setcounter{page}{1}
\setcounter{section}{1}
\renewcommand{\theequation}{A.\arabic{equation}}
\renewcommand{\thesubsection}{B.\arabic{subsection}}
\setcounter{equation}{0}

\begin{center}{\bf \Large Supplementary Material to ``RANK: Large-Scale Inference with Graphical Nonlinear Knockoffs"}

\bigskip

Yingying Fan, Emre Demirkaya, Gaorong Li and Jinchi Lv
\end{center}

\noindent This Supplementary Material contains additional technical details for the proofs of Lemmas \ref{error-lem}--\ref{lem:3-prime}. All the notation is the same as in the main body of the paper.

\section{Additional technical details} \label{appB}

\subsection{Lemma \ref{error-lem} and its proof} \label{appB.1}

\begin{lemma}\label{error-lem}
Assume that $\bX = (X_{ij})\in \mathbb{R}^{n\times p}$ has independent rows with distribution $N(\bzero, \bSig_0)$, $\Lambda_{\max}(\bSig_0) \leq M $, and $\bveps = (\veps_1,\cdots, \veps_n)^T$ has i.i.d. components with $ \mathbb{P} \{ | \veps _{i} | > t \} \leq C_{1} \exp ( -C_{1}^{-1} t^{2} )$ for $t > 0$ and some constants $M, C_{1} > 0$.
Then we have
\begin{eqnarray*}
\mathbb{P}\left\{\Big\|\frac{1}{n}{\bX}^T\bveps\Big\|_{\infty}\le C\sqrt{(\log p)/n}\right\}\geq 1- p^{-c}
\end{eqnarray*}
for some constant $c > 0$ and large enough constant $C > 0$.
\end{lemma}

\noindent \textit{Proof}. First observe that $ \mathbb{P} (| X_{ij}| > t ) \leq 2 \exp \{ - (2 M )^{-1} t^{2} \} $ for $t > 0$,
since $X_{ij} \sim N(0, \bSig_{0,jj}) $ and $\bSig_{0,jj} \leq \Lambda_{\max}(\bSig_0) \leq M$, where $\bSig_{0,jj}$ denotes the $j$th diagonal entry of matrix $\bSig_0$. By assumption, we also have $ \mathbb{P} ( | \veps _{i} | > t ) \leq C_{1} \exp \{ -C_{1}^{-1} t^{2} \}$. Combining these two inequalities yields
\begin{align*}
\mathbb{P} ( | \veps_{i} X_{ij} | > t )
& \leq \mathbb{P}(|\veps_{i}| > \sqrt{t}) + \mathbb{P}(|X_{ij}| > \sqrt{t}) \\
& \leq C_{1} \exp \{ -C_{1}^{-1} t \} + 2 \exp \{ - (2 M )^{-1} t \} \\
& \leq C_{2} \exp\{-C_{2}^{-1} t\},
\end{align*}
where $C_{2} > 0$ is some constant that depends only on constants $C_{1}$ and $M$. Thus by Lemma 6 in \cite{Fan2016interaction}, there exists some constant $\widetilde{C}_{1} > 0$ such that
\begin{align}\label{eq: e066}
\mathbb{P}(|n^{-1} \sum_{i=1}^{n}\veps_{i} X_{ij}| > z )
\leq \widetilde{C}_{1} \exp \{ - \widetilde{C}_{1} n z^{2})
\end{align}
for all $0 < z < 1$.

Denote by ${\bX}_{j}$ the $j$th column of matrix ${\bX}$. Then by \eqref{eq: e066}, the union bound leads to
\begin{align*}
1 - \mathbb{P}\left(\Big\|n^{-1}{\bX}^T\bveps\Big\|_{\infty} \le z\right)
& = \mathbb{P}\left(\Big\|n^{-1}{\bX}^T\bveps\Big\|_{\infty}> z\right) \\
& = \mathbb{P}\Big(\max_{1\le j\le p}|n^{-1}\bveps^T{\bX}_{j}|> z\Big)\\
& \le \sum_{j=1}^p \mathbb{P}(|n^{-1} \sum_{i=1}^{n}\veps_{i} X_{ij}| > z )\\
& \leq p \widetilde{C}_{1} \exp \{ - \widetilde{C}_{1} n z^{2}).
\end{align*}
Letting  $z =  C\sqrt{(\log p)/n} $ in the above inequality, we obtain
\begin{align*}
\mathbb{P} \Big( \Big\|n^{-1}{\bX}^T\bveps\Big\|_{\infty} \le C\sqrt{(\log p)/n}\Big)
\geq 1 - \widetilde{C}_{1} p^{-(\widetilde{C}_{1} C^{2} -1)}.
\end{align*}
Taking large enough positive constant $C$ completes the proof of Lemma \ref{error-lem}.

\subsection{Lemma \ref{lem: beta-est} and its proof} \label{appB.2}

\begin{lemma}\label{lem: beta-est}
Assume that all the conditions of Proposition \ref{pro1} hold and $a_n [(L_p+L_p')^{1/2}+K_n^{1/2}] = o(1)$. Then we have
	\begin{align*}
	P\left\{\sup_{\sbOmg \in \calA, \, |\mathcal{S}|\leq K_n}\left\| \widetilde{\brho}_\mathcal{S}-\widetilde{\bG}_{\mathcal{S},\mathcal{S}}\bbeta_{\mathbb{T},\mathcal{S}}\right\|_\infty \leq C_4\sqrt{(\log p)/n}\right\} = 1 - O(p^{-c_4})
	\end{align*}
for some constants $c_4, C_4>0$.
\end{lemma}

\noindent \textit{Proof}. In this proof, we use $c$ and $C$ to denote generic positive constants and use the same notation as in the proof of Proposition \ref{pro1} in Section \ref{appA.6}. Since $\bbeta_{\mathbb{T}}=(\bbeta_0^T, 0,\ldots,0)^T$ with $\bbeta_0$ the true regression coefficient vector,  it is easy to check that $\wt{\bX}_{\rm KO}\bbeta_{\mathbb{T}}=\bX\bbeta_0$. In view of $\by = \bX\bbeta_0 + \bveps$, it follows from the definitions of $\widetilde{\brho}$  and $\widetilde{\bG}$ that
	\begin{eqnarray*}
		\widetilde{\brho}_\mathcal{S}-\widetilde{\bG}_{\mathcal{S},\mathcal{S}}\bbeta_{\mathbb{T},\mathcal{S}}
		&=&\frac{1}{n}\wt{\bX}_{\rm KO, \mathcal{S}}^T\bX\bbeta_0+\frac{1}{n}\wt{\bX}_{\rm KO, \mathcal{S}}^T\bveps-\dfrac{1}{n}\wt{\bX}_{\rm KO, \mathcal{S}}^{T}\wt{\bX}_{\rm KO, \mathcal{S}}\bbeta_{\mathbb{T},\mathcal{S}}\\
		&=&\frac{1}{n}{\bX}_{\rm KO, \mathcal{S}}^T\bveps+\frac{1}{n}(\wt{\bX}_{\rm KO, \mathcal{S}}-{\bX}_{\rm KO, \mathcal{S}})^T\bveps.
	\end{eqnarray*}
	Using the triangle inequality, we deduce
	\begin{eqnarray*}
		\|\widetilde{\brho}_\mathcal{S}-\widetilde{\bG}_{\mathcal{S},\mathcal{S}}\bbeta_{\mathbb{T},\mathcal{S}}\|_{\infty}\le \Big\|\frac{1}{n}{\bX}_{\rm KO, \mathcal{S}}^T\bveps\Big\|_{\infty}+\Big\|\frac{1}{n}(\wt{\bX}_{\rm KO, \mathcal{S}}-{\bX}_{\rm KO, \mathcal{S}})^T\bveps\Big\|_{\infty}.
	\end{eqnarray*}
We will bound both terms on the right hand side of the above inequality.

By Lemma \ref{error-lem},  we can show that for the first term,
\[ \Big\|\dfrac{1}{n}{\bX}_{\rm KO, \mathcal{S}}^T\bveps\Big\|_{\infty} \le\Big\|\dfrac{1}{n}{\bX}_{\rm KO}^T\bveps\Big\|_{\infty}\le C\sqrt{(\log p)/n} \]
with probability at least $1-p^{-c}$ for some constants $C, c > 0$. We will prove that with probability at least $1-o(p^{-c})$,
	\begin{align}\label{eq: e046}
	\Big\|\dfrac{1}{n}(\wt{\bX}_{\rm KO, \mathcal{S}}-{\bX}_{\rm KO, \mathcal{S}})^T\bveps\Big\|_{\infty} \leq Ca_n(L_p+L_p^{\prime})^{1/2}\sqrt{(\log p)/n} + Ca_n\sqrt{n^{-1}K_n(\log p)}.
	\end{align}
	Then the desired result in this lemma can be shown by noting that $a_n[(L_p+L_p^{\prime})^{1/2}+K_n^{1/2}] \rightarrow 0$.
	
It remains to prove \eqref{eq: e046}. Recall that matrices $\breve{\bX}_{\mathcal S}$ and $\breve\bX_{0, \mathcal{S}}$ can be written as
	\begin{align*}
	&\breve{\bX}_{\mathcal S} = \bX( \bI- \bOmg\diag\{\bs\})_{\mathcal S} +\bZ\bB_{0,\mathcal{S}}(\bB_{0,\mathcal{S}}^T\bB_{0,\mathcal{S}})^{-1/2} \Big(({\bB}_{\mathcal S}^{\sbOmg})^T{\bB}_{\mathcal S}^{\sbOmg}\Big)^{1/2}, \\
	&
	\breve\bX_{0,\mathcal S} = \bX( \bI-\bOmg_0 \diag\{\bs\})_{\mathcal S} + \bZ\bB_{0,\mathcal{S}},
	\end{align*}
	where the notation is the same as in the proof of Proposition \ref{pro1} in Section \ref{appA.6}.
	By the definitions of $\wt{\bX}_{\rm KO}$ and ${\bX}_{\rm KO}$, it holds that
	\begin{align}\label{eq: e044}
	\Big\|\frac{1}{n}(\wt{\bX}_{\rm KO, \mathcal{S}}-{\bX}_{\rm KO, \mathcal{S}})^T\bveps\Big\|_{\infty} = \Big\|\frac{1}{n}(\breve{\bX}_\mathcal{S}-\breve{\bX}_{0,\mathcal{S}})^T\bveps\Big\|_{\infty},
	\end{align}
	where $\breve{\bX}_\mathcal{S}$ and $\breve{\bX}_{0,\mathcal{S}}$ represent the submatrices formed by columns in $\mathcal{S}$. We now turn to analyzing the term $n^{-1}(\breve{\bX}_\mathcal{S}-\breve{\bX}_{0,\mathcal{S}})^T\bveps$.
Some routine calculations give
		\begin{eqnarray*}		\frac{1}{n}(\breve{\bX}_\mathcal{S}-{\breve\bX_{0,\mathcal{S}}})^T\bveps & =&\dfrac{1}{n}\Big(\big((\bOmg_0-{\bOmg}){\rm diag}\{\bs\}\big)_{\mathcal S}\Big)^T\bX^T\bveps \\
			& & +\dfrac{1}{n}\Big( \big(({\bB}_{\mathcal S}^{\sbOmg})^T{\bB}_{\mathcal S}^{\sbOmg}\big)^{1/2}(\bB_{0,\mathcal{S}}^T\bB_{0,\mathcal{S}})^{-1/2}-\bI\Big)\bB_{0,\mathcal{S}}^T\bZ^T\bveps.
		\end{eqnarray*}
Thus it follows from  $s_j\le2\Lambda_{\max}(\bSig_0)$ for all $1 \leq j \leq p$ and the triangle inequality that
		\begin{eqnarray}\label{A-11}
		\nonumber	\Big\|\frac{1}{n}(\breve{\bX}_\mathcal{S}-{\breve\bX_{0,\mathcal{S}}})^T\bveps\Big\|_{\infty}
		&\le& 2\Lambda_{\min}(\bSig_0)\Big\|\dfrac{1}{n}(\bOmg_{0,\mathcal{S}}-{\bOmg}_\mathcal{S})^T\bX^T\bveps\Big\|_{\infty}\\
		&+ &\Big\|\dfrac{1}{n}\Big( \big(({\bB}_{\mathcal S}^{\sbOmg})^T{\bB}_{\mathcal S}^{\sbOmg}\big)^{1/2}(\bB_{0,\mathcal{S}}^T\bB_{0,\mathcal{S}})^{-1/2}-\bI\Big)\bB_{0,\mathcal{S}}^T\bZ^T\bveps\Big\|_{\infty}.
		\end{eqnarray}
	
We first examine the upper bound for  $\Big\|\dfrac{1}{n}(\bOmg_{0,\mathcal{S}}-{\bOmg}_\mathcal{S})^T\bX^T\bveps\Big\|_{\infty}$  in \eqref{A-11}. Since $\bOmg \in \calA$ and $\bOmg_0$ is $L_p$-sparse, by Lemma \ref{error-lem} we deduce
		\begin{eqnarray}
		\nonumber		\Big\|\dfrac{1}{n}(\bOmg_{0,\mathcal{S}}-{\bOmg}_\mathcal{S})^T\bX^T\bveps\Big\|_{\infty}&\le&\Big\|\dfrac{1}{n}(\bOmg_0-{\bOmg})\bX^T\bveps\Big\|_{\infty}\\
		&\le&\|\bOmg_0-{\bOmg}\|_1\Big\|\dfrac{1}{n}\bX^T\bveps\Big\|_{\infty} \nonumber\\
		&\le&\sqrt{L_p+L_p^{\prime}}\|\bOmg - \bOmg_0\|_{2}\cdot C\sqrt{(\log p)/n}\nonumber\\
		&\le&Ca_n(L_p+L_p^{\prime})^{1/2}\sqrt{(\log p)/n}.\label{eq: e045}
		\end{eqnarray}
We can also bound the second term on the right hand side of \eqref{A-11} as
		\begin{align*}
		& \Big\|\dfrac{1}{n}\Big( \big(({\bB}_{\mathcal S}^{\sbOmg})^T{\bB}_{\mathcal S}^{\sbOmg}\big)^{1/2}(\bB_{0,\mathcal{S}}^T\bB_{0,\mathcal{S}})^{-1/2}-\bI\Big)\bB_{0,\mathcal{S}}^T\bZ^T\bveps\Big\|_{\infty}\\
		&\leq 	\Big\|\big(({\bB}_{\mathcal S}^{\sbOmg})^T{\bB}_{\mathcal S}^{\sbOmg}\big)^{1/2}(\bB_{0,\mathcal{S}}^T\bB_{0,\mathcal{S}})^{-1/2}-\bI\Big\|_1\Big\|\frac{1}{n}\bB_{0,\mathcal{S}}^T\bZ^T\bveps\Big\|_{\infty}\\
		&\leq \sqrt{2|\mathcal S|}\Big\|\big(({\bB}_{\mathcal S}^{\sbOmg})^T{\bB}_{\mathcal S}^{\sbOmg}\big)^{1/2}(\bB_{0,\mathcal{S}}^T\bB_{0,\mathcal{S}})^{-1/2}-\bI\Big\|_2 \Big\|\frac{1}{n}\bB_{0,\mathcal{S}}^T\bZ^T\bveps\Big\|_{\infty}\\
		& \leq \sqrt{2K_n}Ca_n \sqrt{(\log p)/n} = Ca_n\sqrt{n^{-1}K_n(\log p)},
		\end{align*}
		where the second to the last step is entailed by Lemma \ref{lemma: 1} in Section \ref{appA.3}  and Lemma \ref{lem: 3} in Section \ref{appB.3}. Therefore, combining this inequality with \eqref{eq: e044}--\eqref{eq: e045} results in  \eqref{eq: e046}, which concludes the proof of Lemma \ref{lem: beta-est}.

\subsection{Lemma \ref{lem: 3} and its proof} \label{appB.3}

\begin{lemma}\label{lem: 3}
	Under the conditions of Proposition \ref{pro1}, it holds that with probability at least $1-O(p^{-c})$,
	\[
	\sup_{|\mathcal S| \leq K_n}	\Big\|\frac{1}{n}\bB_{0,\mathcal{S}}^T\bZ^T\bveps\Big\|_{\infty} \geq C\sqrt{(\log p)/n}
	\]
for some constant $C > 0$.
\end{lemma}

\noindent \textit{Proof}. Since this is a specific case of Lemma \ref{lem:3-prime} in Section \ref{appB.6}, the proof is omitted.

\subsection{Lemma \ref{lem: cond-Shat} and its proof} \label{appB.4}

\begin{lemma}\label{lem: cond-Shat}
Under the conditions of Proposition \ref{pro1} and Theorem \ref{lem: cond-Shat0}, there exists some constant $c \in (2(qs)^{-1},1)$ such that with asymptotic probability one, $|\widehat{\mathcal{S}}^{\sbOmg}| \geq cs$
holds uniformly over all $\bOmg \in \calA$ and $|\mathcal{S}|\leq K_n$, where $\widehat{\mathcal{S}}^{\sbOmg} = \{j:  W_j^{\sbOmg, \mathcal{S}}\geqslant T\}$.
\end{lemma}

\noindent \textit{Proof}. Again we use $c$ and $C$ to denote generic positive constants whose values may change from line to line.  By Proposition \ref{pro1} in Section \ref{appA.6}, we have with probability at least  $1-O(p^{-c_1})$ that uniformly over all $\bOmg\in \calA$ and $|\mathcal S| \leq K_n$,
	\begin{align*}
	\max_{1\leq j\leq p}|\hbeta_{j}(\lambda; {\bOmg}, \mathcal S) - \beta_{0,j}| \leq C\sqrt{sn^{-1}(\log p)} \text{ and } \max_{1 \leq j \leq  p}|\hbeta_{j+p}(\lambda; {\bOmg}, \mathcal S)|\leq C\sqrt{sn^{-1}(\log p)}
	\end{align*}
for some constants $C, c_1>0$. Thus for each $1 \leq j \leq p$, we have
\begin{align}\label{eq: e60prime}
\nonumber W_j^{\sbOmg, \mathcal S} & = |\hbeta_{j}(\lambda; {\bOmg}, \mathcal S) | - |\hbeta_{j+p}(\lambda; {\bOmg}, \mathcal S)| \\
& \geq - |\hbeta_{j+p}(\lambda; {\bOmg}, \mathcal S)| \geq - C \sqrt{s n^{-1} (\log p)}.
	\end{align}

On the other hand, for each $j \in \mathcal{S}_{2} = \{j: \beta_{0,j}\gg \sqrt{sn^{-1}(\log p)}\}$ it holds that
	\begin{align}\label{eq: e060}
\nonumber W_j^{\sbOmg, \mathcal S} & = |\hbeta_{j}(\lambda; {\bOmg}, \mathcal S) | - |\hbeta_{j+p}(\lambda; {\bOmg}, \mathcal S)| \\
& \geq   |\beta_{0,j}| - |\hbeta_{j}(\lambda; {\bOmg}, \mathcal S) - \beta_{0,j}| - |\hbeta_{j+p}(\lambda; {\bOmg}, \mathcal S)| \gg C \sqrt{s n^{-1} (\log p)}.
	\end{align}
Thus in order for  any $W_j^{\sbOmg, \mathcal S}$, $1 \leq j \leq p$ to fall below $-T$, we must have $W_j^{\sbOmg, \mathcal S} \geq T$ for all $j\in \mathcal{S}_2$. This entails that
\begin{align} \label{eqn:TzerotoStwo}
\left|\{j: W_j^{\sbOmg, \mathcal S} \geq T\}\right| \geq |\mathcal{S}_2| \geq cs,
\end{align}
 which completes the proof of Lemma \ref{lem: cond-Shat}.

\subsection{Lemma \ref{bound-lem} and its proof} \label{appB.5}

\begin{lemma}\label{bound-lem}
Assume that all the conditions of Proposition \ref{pro1} hold and $a_{2n} = a_n +( L_{p}' + K_n) \{(\log p)/n\}^{1/2} = o(1)$. Then it holds that
\begin{align*}
P\left\{\sup_{\sbOmg \in \calA, \, |\mathcal{S}|\leq K_n}\left\| \widetilde{\bG}_{\mathcal{S},\mathcal{S}} - \bG_{\mathcal{S},\mathcal{S}} \right\|_{\max} \leq C_8 a_{2,n} \right\} = 1 - O(p^{-c_8})
\end{align*}
for some constants $c_8, C_8>0$.
\end{lemma}

\noindent \textit{Proof}. In this proof, we adopt the same notation as used in the proof of Proposition \ref{pro1} in Section \ref{appA.6}. In light of  (\ref{eqn:def-Gtilde}), we have $\widetilde{\bG} = n^{-1} [\bX, \breve{\bX}^{\sbOmg}]^{T}[\bX, \breve{\bX}^{\sbOmg}]$. Thus
the matrix difference $\widetilde{\bG}_{\mathcal{S},\mathcal{S}} - \bG_{\mathcal{S},\mathcal{S}} $ can be represented in block form as
\begin{align*}
\widetilde{\bG}_{\mathcal{S},\mathcal{S}} - \bG_{\mathcal{S},\mathcal{S}}
& =
\dfrac{1}{n}
\left(\begin{array}{cc}
\bX^{T}_{\mathcal S} \bX_{\mathcal S} & (\breve{\bX}_{\mathcal S}^{\sbOmg})^{T} \bX_{\mathcal S}  \\
\bX_{\mathcal S}^{T}\breve{\bX}_{\mathcal S}^{\sbOmg} & (\breve{\bX}_{\mathcal S}^{\sbOmg})^{T} \breve{\bX}_{\mathcal S}^{\sbOmg}
\end{array}\right)
-
\left(\begin{array}{cc}
\bSig_0 & \bSig_0 - \diag\{\bs\} \\
\bSig_0 - \diag\{\bs\} & \bSig_0
\end{array} \right)
_{\mathcal{S},\mathcal{S}}
\\
& =
\left(\begin{array}{cc}
n^{-1} \bX^{T}_{\mathcal S} \bX_{\mathcal S}
- \bSig_{0, \mathcal{S}, \mathcal{S}}
& n^{-1} (\breve{\bX}_{\mathcal S}^{\sbOmg})^{T} \bX_{\mathcal S}
- \big(\bSig_0 - \diag\{\bs\}\big)_{\mathcal{S},\mathcal{S}}
\\
n^{-1} \bX_{\mathcal S}^{T}\breve{\bX}_{\mathcal S}^{\sbOmg}
- \big(\bSig_0 - \diag\{\bs\}\big)_{\mathcal{S},\mathcal{S}}
& n^{-1} (\breve{\bX}_{\mathcal S}^{\sbOmg})^{T} \breve{\bX}_{\mathcal S}^{\sbOmg} - \bSig_{0, \mathcal{S}, \mathcal{S}}
\end{array}\right).
\end{align*}
Note that the off-diagonal blocks are the transposes of each other. Then we see that
$ \| \widetilde{\bG}_{\mathcal{S}, \mathcal{S}} - \bG_{\mathcal{S}, \mathcal{S}} \|_{\max} $ can be bounded by the maximum of $\|\eta_1\|_{\max}$, $\|\eta_2\|_{\max}$, and $\|\eta_3\|_{\max}$ with
\begin{align*}
\eta_1 & = n^{-1} \bX^{T}_{\mathcal S} \bX_{\mathcal S}
- \bSig_{0, \mathcal{S}, \mathcal{S}} , \\
\eta_2 & = n^{-1} \bX_{\mathcal S}^{T}\breve{\bX}_{\mathcal S}^{\sbOmg}
- \big(\bSig_0 - \diag\{\bs\}\big)_{\mathcal{S},\mathcal{S}} , \\
\eta_3 & = n^{-1} (\breve{\bX}_{\mathcal S}^{\sbOmg})^{T} \breve{\bX}_{\mathcal S}^{\sbOmg} - \bSig_{0, \mathcal{S}, \mathcal{S}}.
\end{align*}
To bound these three terms, we define three events
\begin{align*}
\mathcal{E}_5 & =
\Big\{ \| n^{-1} \bX^{T}\bX - \bSig_{0} \|_{\max}
\leq C \sqrt{(\log p)/n} \Big\}, \\
\mathcal{E}_6 & =
\Big\{ \sup_{|\mathcal S| \leq K_n} \Big\| n^{-1} \bB_{0,\mathcal{S}}^T\bZ^T\bX\Big\|_{\infty} \leq C\sqrt{(\log p)/n} \Big\}, \\
\mathcal{E}_7 & =
\Big\{ \sup_{|\mathcal S| \leq K_n}
\Big\| n^{-1} \bB_{0,\mathcal{S}}^T\bZ^T \bZ \bB_{0,\mathcal{S}} - \bB_{0,\mathcal{S}}^T \bB_{0,\mathcal{S}} \Big\|_{\max}
\leq C \sqrt{(\log p)/n} \Big\}.
\end{align*}
By Lemma \ref{lem:3-prime} in Section \ref{appB.6}, it holds that $ P( \mathcal{E}_6 ) \geq 1 - O(p^{-c}) $ and $ P( \mathcal{E}_7 ) \geq 1 - O(p^{-c}) $.
Using Lemma A.3 in \cite{Bickel2008}, we also have $ P( \mathcal{E}_5 ) \geq 1 - O(p^{-c}) $.
Combining these results yields
\[ P( \mathcal{E}_5 \cap \mathcal{E}_6 \cap \mathcal{E}_7) \geq 1 - O(p^{-c}) \]
with $c > 0$ some constant.

Let us first consider term $\eta_{1}$. Conditional on $\mathcal{E}_5$, it is easy to see that
\begin{align} \label{eqn:eta1}
\| \eta_{1} \|_{\max} \leq \| n^{-1} \bX^{T}\bX - \bSig_{0} \|_{\max}
\leq C \sqrt{(\log p)/n}.
\end{align}
We next bound $\| \eta_{2} \|_{\max} $ conditional on $\mathcal{E}_5 \cap \mathcal{E}_6$.
To simplify the notation, denote by $\widetilde{\bB}^{\mathcal{S}, \sbOmg} = (\bB_{0,\mathcal{S}}^T\bB_{0,\mathcal{S}})^{-1/2} \Big(({\bB}_{\mathcal S}^{\sbOmg})^T{\bB}_{\mathcal S}^{\sbOmg}\Big)^{1/2}
$. By the definition of $\breve{\bX}_{\mathcal S}$,
we deduce
\begin{align*}
\eta_{2}
& =
n^{-1} \bX_{\mathcal S}^{T}\breve{\bX}_{\mathcal S}^{\sbOmg}
- \big(\bSig_0 - \diag\{\bs\}\big)_{\mathcal{S},\mathcal{S}}
\\
& =
n^{-1} \bX_{\mathcal S}^{T} \bX( \bI- \bOmg\diag\{\bs\})_{\mathcal S} + n^{-1} \bX_{\mathcal S}^{T} \bZ\bB_{0,\mathcal{S}}\widetilde{\bB}^{\mathcal{S}, \sbOmg}
- \big(\bSig_0 - \diag\{\bs\}\big)_{\mathcal{S},\mathcal{S}}
\\
& =
\big( ( n^{-1} \bX^{T} \bX - \bSig_0 )
( \bI - \bOmg\diag\{\bs\})
\big)_{\mathcal{S},\mathcal{S}}
+
\big( \diag\{\bs\} - \bSig_0 \bOmg\diag\{\bs\} \big)_{\mathcal{S},\mathcal{S}}
+
n^{-1} \bX_{\mathcal S}^{T} \bZ\bB_{0,\mathcal{S}}\widetilde{\bB}^{\mathcal{S}, \sbOmg}
\\
& \equiv \eta_{2,1} + \eta_{2,2} + \eta_{2,3}.
\end{align*}
We will examine the above three terms separately.

Since $\bOmg$ is $L_p'$-sparse, $\|\bI-\bOmg_{0}\diag(\bs)\|_2\leq \|\bI\|_2 + \|\bOmg_{0}\diag(\bs)\|_2 \leq C$, and $\|(\bOmg - \bOmg_{0}) \diag\{\bs\}\|_2\leq Ca_{2n}$, we have
\begin{align}
\Big\| \bI - \bOmg\diag\{\bs\} \Big\|_{1}
& \le \sqrt{L_{p}'} \Big\| \bI - \bOmg\diag\{\bs\} \Big\|_{2}
\nonumber \\
& \le \sqrt{L_{p}'} \Big(\big\| \bI - \bOmg_{0} \diag\{\bs\} \big\|_{2} +
\big\| (\bOmg - \bOmg_{0}) \diag\{\bs\} \big\|_{2} \Big)
\nonumber \\
& \le C \sqrt{L_{p}'}. \label{eqn:Eye-OmgS}
\end{align}
Thus it follow from (\ref{eqn:Eye-OmgS}) that conditional on $\mathcal{
E}_{5}$,
\begin{align}
\| \eta_{2,1} \|_{\max} & =
\Big\| \big( ( n^{-1} \bX^{T} \bX - \bSig_0 )
( \bI - \bOmg\diag\{\bs\})
\big)_{\mathcal{S},\mathcal{S}} \Big\|_{\max}
\nonumber \\
& \le
\Big\|  ( n^{-1} \bX^{T} \bX - \bSig_0 ) ( \bI - \bOmg\diag\{\bs\}) \Big\|_{\max}
\nonumber \\
& \le
\Big\|  n^{-1} \bX^{T} \bX - \bSig_0  \Big\|_{\max} \Big\| \bI - \bOmg\diag\{\bs\} \Big\|_{1}
\nonumber \\
& \le  C \sqrt{L_{p}'}\sqrt{(\log p)/n}. \label{eqn:eta2-1}
\end{align}
For term $\eta_{2, 2}$, it holds that
\begin{align}
\| \eta_{2,2} \|_{\max} & =
\Big\| \big( \diag\{\bs\} - \bSig_0 \bOmg\diag\{\bs\} \big)_{\mathcal{S},\mathcal{S}} \Big\|_{\max}
\nonumber \\
& \le
C \|  \bI - \bSig_0 \bOmg  \|_{\max}
\le C \| \bSig_0 \|_{2} \|  \bOmg_0 - \bOmg  \|_{2}
\le C a_{n}. \label{eqn:eta2-2}
\end{align}

Note that by Lemma \ref{lemma: 1} in Section \ref{appA.3}, we have
\[ \| \widetilde{\bB}^{\mathcal{S}, \sbOmg} \|_{1}
\le \sqrt{|\mathcal{S}|} \| \widetilde{\bB}^{\mathcal{S}, \sbOmg} \|_{2}
\le \sqrt{|\mathcal{S}|} ( \| \widetilde{\bB}^{\mathcal{S}, \sbOmg} - \bI \|_{2} + 1 )
\le C \sqrt{|\mathcal{S}|}
\le C \sqrt{K_n} \]
when $|\mathcal S| \leq K_n$. Then conditional on $\mathcal{E}_{6}$, it holds that
\begin{align}
\| \eta_{2,3} \|_{\max} & =
\| n^{-1} \bX_{\mathcal S}^{T} \bZ \bB_{0,\mathcal{S}}\widetilde{\bB}^{\mathcal{S}, \sbOmg} \|_{\max}
\nonumber \\
& \le
\| n^{-1} \bX_{\mathcal S}^{T} \bZ \bB_{0,\mathcal{S}} \|_{\max} \| \widetilde{\bB}^{\mathcal{S}, \sbOmg} \|_{1}
\nonumber \\
& \le C  \sqrt{n^{-1}K_n(\log p)}. \label{eqn:eta2-3}
\end{align}
Thus combining (\ref{eqn:eta2-1})--(\ref{eqn:eta2-3}) leads to
\begin{equation} \label{eqn:eta2}
\| \eta_2 \|_{\max} \le C \{a_n +  \sqrt{n^{-1}L_{p}'(\log p)} +  \sqrt{n^{-1}K_n(\log p)}\}.
\end{equation}

We finally deal with term $\eta_{3}$. Some routine calculations show that
\begin{align*}
\eta_{3}
& =
n^{-1} (\breve{\bX}_{\mathcal S}^{\sbOmg})^{T} \breve{\bX}_{\mathcal S}^{\sbOmg} - \bSig_{0, \mathcal{S}, \mathcal{S}}.
\\
& =
n^{-1}
\big( ( \bI- \bOmg\diag\{\bs\})_{\mathcal S}^{T} \bX^{T}  + ( \widetilde{\bB}^{\mathcal{S}, \sbOmg} )^{T} \bB_{0,\mathcal{S}}^{T} \bZ^{T}
\big) \big( \bX( \bI- \bOmg\diag\{\bs\})_{\mathcal S} + \bZ\bB_{0,\mathcal{S}}\widetilde{\bB}^{\mathcal{S}, \sbOmg}
\big)
- \bSig_{0, \mathcal{S}, \mathcal{S}}
\\
& =
\Big( n^{-1}
(\bI- \bOmg\diag\{\bs\})^{T} \bX^{T}
\bX( \bI- \bOmg\diag\{\bs\})
- \bSig_{0} + \bB_{0}^T \bB_{0}  \Big)_{\mathcal{S}, \mathcal{S}}
\\
& + n^{-1}
( \widetilde{\bB}^{\mathcal{S}, \sbOmg} )^{T} \bB_{0,\mathcal{S}}^{T} \bZ^{T}
\bX( \bI- \bOmg\diag\{\bs\})_{\mathcal S}
+
( \bI- \bOmg\diag\{\bs\})_{\mathcal S}^{T} \bX^{T}
\bZ\bB_{0,\mathcal{S}}\widetilde{\bB}^{\mathcal{S}, \sbOmg}
\\ & +
\big( ( \widetilde{\bB}^{\mathcal{S}, \sbOmg} )^{T} \bB_{0,\mathcal{S}}^{T} \bZ^{T}
\bZ\bB_{0,\mathcal{S}}\widetilde{\bB}^{\mathcal{S}, \sbOmg} - \bB_{0,\mathcal{S}}^T \bB_{0,\mathcal{S}} \big)
\\
& \equiv \eta_{3,1} + \eta_{3,2} + \eta_{3,2}^{T} + \eta_{3,3}.
\end{align*}
Conditional on event $\mathcal{E}_{5}$, with some simple matrix algebra we derive
\begin{align}
\| \eta_{3,1}\| & = \Big\| \Big( n^{-1}
(\bI- \bOmg\diag\{\bs\})^{T} \bX^{T}
\bX( \bI- \bOmg\diag\{\bs\})
- \bSig_{0} + \bB_{0}^T \bB_{0}  \Big)_{\mathcal{S}, \mathcal{S}} \Big\|_{\max}
\nonumber \\
& \le
\Big\| n^{-1}
(\bI- \bOmg\diag\{\bs\})^{T} \bX^{T}
\bX( \bI- \bOmg\diag\{\bs\})
- \bSig_{0} + \bB_{0}^T \bB_{0} \Big\|_{\max}
\nonumber \\
& \le
\Big\|
(\bI- \bOmg\diag\{\bs\})^{T} (n^{-1}\bX^{T}
\bX - \bSig_{0}) (\bI- \bOmg\diag\{\bs\})
\Big\|_{\max}
\nonumber \\ & +
\Big\|
(\bI- \bOmg\diag\{\bs\})^{T} \bSig_{0} (\bI- \bOmg\diag\{\bs\})
- \bSig_{0}
+ 2 \diag\{\bs\}
- \diag\{\bs\}\bOmg_{0}\diag\{\bs\}
\Big\|_{\max}
\nonumber \\
& \le
\|
n^{-1}\bX^{T} \bX - \bSig_{0} \|_{\max}
\| (\bI- \bOmg\diag\{\bs\}) \|^{2}_{1}
\nonumber \\ & +
\| \diag\{\bs\} (\bI- \bOmg \Sig_0) \|_{\max}
+\| (\bI- \bOmg \Sig_0) \diag\{\bs\} \|_{\max}
+\| \diag\{\bs\} (\bOmg_0 - \bOmg \Sig_0 \bOmg) \diag\{\bs\} \|_{\max}
\nonumber \\
& \leq
C L_{p}' \sqrt{(\log p)/n} + C a_{n}, \label{eqn:eta3-1}
\end{align}
where the last step used (\ref{eqn:Eye-OmgS}) and calculations similar to (\ref{eqn:eta2-2}).

It follows from (\ref{eqn:Eye-OmgS}) and the previously proved result $\| \widetilde{\bB}^{\mathcal{S}, \sbOmg} \|_{1}
\le C \sqrt{K_n}$ for $|\mathcal{S}| \le K_n$ that conditional on event $\mathcal{E}_{6}$,
\begin{align}
\| \eta_{3,2}\| & = \| n^{-1} ( \widetilde{\bB}^{\mathcal{S}, \sbOmg} )^{T} \bB_{0,\mathcal{S}}^{T} \bZ^{T} \bX( \bI- \bOmg\diag\{\bs\})_{\mathcal S} \|_{\max}
\nonumber \\
& \leq
\|  \widetilde{\bB}^{\mathcal{S}, \sbOmg} \|_{1}
\| n^{-1} \bB_{0,\mathcal{S}}^{T} \bZ^{T} \bX \|_{\max}
\| ( \bI- \bOmg\diag\{\bs\})_{\mathcal S} \|_{1}
\nonumber \\
& \leq
C \sqrt{K_n} \sqrt{ L_p'n^{-1} (\log p)}
\nonumber \\
& =
C  \sqrt{n^{-1} K_n L_p'(\log p)}. \label{eqn:eta3-2}
\end{align}
Finally, by Lemma \ref{lemma: 1} it holds that conditioned on $\mathcal{E}_{7}$,
\begin{align}
\| \eta_{3,3}\| & =  \Big\| n^{-1}
( \widetilde{\bB}^{\mathcal{S}, \sbOmg} )^{T} \bB_{0,\mathcal{S}}^{T} \bZ^{T} \bZ\bB_{0,\mathcal{S}}\widetilde{\bB}^{\mathcal{S}, \sbOmg} - \bB_{0,\mathcal{S}}^T \bB_{0,\mathcal{S}} \Big\|_{\max}
\nonumber \\
& \leq
\Big\|
( \widetilde{\bB}^{\mathcal{S}, \sbOmg} )^{T}
(n^{-1} \bB_{0,\mathcal{S}}^{T} \bZ^{T} \bZ \bB_{0,\mathcal{S}} - \bB_{0,\mathcal{S}}^T \bB_{0,\mathcal{S}}) \widetilde{\bB}^{\mathcal{S}, \sbOmg} \Big\|_{\max}
\nonumber \\ & +
\Big\|
( \widetilde{\bB}^{\mathcal{S}, \sbOmg} )^{T} \bB_{0,\mathcal{S}}^{T} \bB_{0,\mathcal{S}}\widetilde{\bB}^{\mathcal{S}, \sbOmg} - \bB_{0,\mathcal{S}}^T \bB_{0,\mathcal{S}} \Big\|_{\max}
\nonumber \\
& \leq
\Big\|
n^{-1} \bB_{0,\mathcal{S}}^{T} \bZ^{T} \bZ \bB_{0,\mathcal{S}} - \bB_{0,\mathcal{S}}^T \bB_{0,\mathcal{S}} \Big\|_{\max}
\| \widetilde{\bB}^{\mathcal{S}, \sbOmg} \|_{1}^{2}
+ C a_{n}
\nonumber \\
& \leq
C  K_n\sqrt{(\log p)/n}
+ C a_{n}. \label{eqn:eta3-3}
\end{align}
Therefore, combining (\ref{eqn:eta3-1})--(\ref{eqn:eta3-3}) results in
\begin{align*}
\| \eta_{3} \|_{\max} & \leq Ca_n + C(L_{p}' + K_n + \sqrt{K_nL_p'}) \sqrt{(\log p)/n} \\
& \leq  Ca_n + 2C(L_{p}' + K_n) \sqrt{(\log p)/n},
\end{align*}
which together with (\ref{eqn:eta1}) and (\ref{eqn:eta2}) concludes the proof of Lemma \ref{bound-lem}.

\subsection{Lemma \ref{lem:3-prime} and its proof} \label{appB.6}

\begin{lemma}\label{lem:3-prime}
	Under the conditions of Proposition \ref{pro1}, it holds that with probability at least $1-O(p^{-c})$,
	\begin{align*}
	&\sup_{|\mathcal S| \leq K_n}	\Big\|\frac{1}{n}\bB_{0,\mathcal{S}}^T\bZ^T \bX \Big\|_{\max} \geq C\sqrt{(\log p)/n},\\
	&\sup_{|\mathcal S| \leq K_n}
	\Big\| n^{-1} \bB_{0,\mathcal{S}}^T\bZ^T \bZ \bB_{0,\mathcal{S}} - \bB_{0,\mathcal{S}}^T \bB_{0,\mathcal{S}}  \Big\|_{\max} \geq C\sqrt{(\log p)/n}
	\end{align*}
for some constants $c, C > 0$.
\end{lemma}

\noindent \textit{Proof}. We still use $c$ and $C$ to denote generic positive constants. We start with proving the first inequality. Observe that
\[
	\sup_{|\mathcal S| \leq K_n}	\Big\|\frac{1}{n}\bB_{0,\mathcal{S}}^T\bZ^T \bX \Big\|_{\max} \leq 	\Big\|\frac{1}{n}\bB_{0}^T\bZ^T \bX \Big\|_{\max}.
\]
Thus it remains to prove
\begin{align}\label{eq: e065}
P\left(	\Big\|\frac{1}{n}\bB_{0}^T\bZ^T \bX \Big\|_{\max} \geq C\sqrt{(\log p)/n}\right) \leq o(p^{-c}).
\end{align}

Let $\bU = \bZ\bB_{0} \in \mathbb{R}^{n\times p}$ and denote by $\bU_j$ the $j$th column of matrix $\bU$. We see that the components of  $\bU_j$ are i.i.d. Gaussian with mean zero and variance $\be_j^T\bB_{0}^T\bB_{0}\be_j$, and the vectors $\bU_j$  are independent of $\bveps$. Let $\wt{\bU}_j = (\be_j^T\bB_{0}^T\bB_{0}\be_j)^{-1/2}\bU_j$.  Then it holds that $\wt{\bU}_j \sim N(\bzero, \bI_{n})$.
Since $X_{ij} \sim N(0, \bSig_{0,jj}) $ and $\bSig_{0,jj} \leq \Lambda_{\max}(\bSig_0) \leq C$ with $C > 0$ some constant, it follows from Bernstein's inequality that for $t>0$,
	\begin{align*}
	\mb{P}\left(	\Big\|\frac{1}{n}\bB_{0}^T\bZ^T\bX \Big\|_{\max} \geq t \|\bB_{0}^T\bB_{0}\|_2 \right) &\leq  \sum_{j=1}^{p}\mb{P}\left(\frac{1}{n}\Big|(\bU_j)^T\bX_{i}\Big| \geq t \|\bB_{0}^T\bB_{0}\|_2 \right)\\
	& \leq  \sum_{i=1}^{p}\mb{P}\left(\frac{1}{n}\Big|(\wt\bU_j)^T\bX_{i}\Big| \geq t \right)\\
	& \leq C p \exp(-Cnt^2).
	\end{align*}
	Taking $t = C\sqrt{(\log p)/n}$  with large enough constant $C>0$ in the above inequality yields  	
	\[
	\mb{P}\left(	\Big\|\frac{1}{n}\bB_{0}^T\bZ^T\bX \Big\|_{\max} \geq C\sqrt{(\log p)/n} \|\bB_{0}^T\bB_{0}\|_2 \right) \leq Cp^{-c}
	\]
for some constant $c > 0$. Thus with probability at least $1-O(p^{-c})$, it holds that
	\begin{align*}
	& \Big\|\frac{1}{n}\bB_{0}^T\bZ^T\bX \Big\|_{\max} \leq C\sqrt{(\log p)/n} \|\bB_{0}^T\bB_{0}\|_2\\
	&  = C\sqrt{(\log p)/n} \|\diag(\bs) -\diag(\bs)\bOmg_{0}\diag(\bs) \|_2\\
	& \leq C\sqrt{(\log p)/n},
	\end{align*}
which establishes \eqref{eq: e065} and thus concludes the proof for the first result.

The second inequality follows from
\begin{align*}
 \sup_{|\mathcal S| \leq K_n}
\Big\| n^{-1} \bB_{0,\mathcal{S}}^T\bZ^T \bZ \bB_{0,\mathcal{S}} - \bB_{0,\mathcal{S}}^T \bB_{0,\mathcal{S}}  \Big\|_{\max} \leq
\Big\| n^{-1} \bB_{0}^T\bZ^T \bZ \bB_{0} - \bB_{0}^T \bB_{0}  \Big\|_{\max}
\end{align*}
and  Lemma A.3 in \cite{Bickel2008}, which completes the proof of Lemma \ref{lem:3-prime}.
\end{document}